\documentclass{elsarticle}
\usepackage{latexsym,amssymb,amsmath}
\usepackage{color}
\usepackage[english]{babel}
\usepackage{amsfonts}

\newtheorem{thm}{Theorem}
\newtheorem{lemma}[thm]{Lemma}
\newdefinition{rmk}{Remark}
\newproof{pf}{Proof}
\newtheorem{definition}[thm]{Definition}

\begin{document}
\title{ The regularized 3D Boussinesq equations with fractional Laplacian and no diffusion}
\author{H. Bessaih}
 \ead{bessaih@uwyo.edu}
 \address{University of Wyoming, Department of Mathematics, Dept. 3036, 1000
East University Avenue, Laramie WY 82071, United States}

 \author{B. Ferrario\corref{cor1}}
 \ead{benedetta.ferrario@unipv.it}
 \address{Universit\`a di Pavia, Dipartimento di Matematica, via
  Ferrata 5, 27100 Pavia, Italy}
\cortext[cor1]{Corresponding author}

\begin{abstract}
In this paper, we study the 3D regularized Boussinesq equations. 
The velocity equation is regularized \`a la Leray 
through a smoothing kernel of order $\alpha$ in the nonlinear term and a 
$\beta$-fractional Laplacian; we consider the critical case $\alpha+\beta=\frac{5}{4}$ 
and we assume  $\frac 12 <\beta<\frac 54$.
The temperature equation is a pure transport equation, where the
transport velocity is regularized through the
same smoothing kernel of order $\alpha$.  
We prove global well posedness when the initial velocity is in $H^r$ 
and the initial temperature is in $H^{r-\beta}$ for 
$r>\max(2\beta,\beta+1)$.
This regularity is enough to prove uniqueness of solutions. 
We also prove a continuous dependence of the solutions on the initial conditions. 
\end{abstract}

\date{\today}
\maketitle

\section{Introduction}
We consider the  Boussinesq system 
in a $d$-dimensional space:
\begin{equation}\label{syst0}
\begin{cases}
\partial_t v+(v \cdot \nabla) v -\nu \Delta v +\nabla p =\theta e_d \\
\partial_t \theta+v\cdot \nabla \theta=0\\
\nabla \cdot v=0
\end{cases}
\end{equation}
where $v=v(t,x)$ denotes the velocity vector field, $p=p(t,x)$  the scalar pressure 
and $\theta=\theta(t,x)$ a scalar quantity, which can represent
either the temperature of the fluid or the concentration of a chemical component; 
$e_d$ is the unit vector $(0,\ldots,0,1)$,
the viscosity $\nu$ is a  positive constant.
Suitable initial conditions $v_0, \theta_0$  
and boundary conditions (if needed) are given.
 
For $d=2$,
the well posedness of system \eqref{syst0} in the plane 
has been studied by several authors under different assumptions on the initial data
(see \cite{HL05,C06,AH07,HK07,DP08,DP09}). For $d=3$, 
very little is known; 
it has been proven that there exists a local smooth
solution. Some regularity criterions to get a global (in time) solution
have been obtained  in \cite{QDY10,GF12}. Otherwise, in the particular case of 
axisymmetric initial data, 
\cite{AHK11} shows the  global well posedness
for the Boussinesq system in the whole space.

To overcome the difficulties of the three-dimensional case, 
different models have been proposed.
For instance, one can regularize the equation for the velocity by putting a 
fractional power of the Laplacian; this 
hyper-dissipative Boussinesq system takes  the form

\begin{equation}\label{syst:alpha}
\begin{cases}
\partial_t v+(v \cdot \nabla) v +\nu(-\Delta)^{\beta}v +\nabla p =\theta e_3 \\
\partial_t \theta+v\cdot \nabla \theta=0\\
\nabla \cdot v=0
\end{cases}
\end{equation}
For $\beta>\frac 54$, \cite{ZW13} proved the global well posedness. 
This result has been improved by Ye \cite{Ye}, allowing $\beta = \frac 54$.

Notice that for zero initial temperature $\theta_0$, 
the Boussinesq  system reduces 
to the Navier-Stokes equations. It is well known that the three-dimensional 
Navier-Stokes equations have either a unique local smooth solution or a global 
weak solution. The questions related to the local 
smooth solution being  global or  the global weak solution being  unique are very challenging problems that are still open since the seminal work of Leray.
For this reason, modifications of different 
types have been considered for the three-dimensional Navier-Stokes equations. 
On one side there is the  hyper-viscous model, i.e.
\eqref{syst:alpha} with zero initial temperature; 
when $\beta \ge \frac 54$, uniqueness of the weak solutions has been proved 
in \cite{Lions} (see Remark 6.11 of Chapter 1) and \cite{MS} .
On the other hand,  Olson and Titi 
in \cite{ot} suggested to regularize the equations by modifying
two terms. For a  particular model
of fluid dynamics, they replaced the dissipative term  
 by  a fractional power of the Laplacian and they regularized the bilinear term 
of vorticity stretching \`a la 
 Leray. The well posedness of those equations
is obtained by asking a balance between the modification of the 
nonlinearity and of the  viscous dissipation; at least 
one of them has to be strong enough, while 
the other might be weak.
Similarly, Barbato, Morandin and Romito in \cite{BMR} considered 
the  Leray-$\alpha$ Navier-Stokes equations with fractional dissipation

\begin{equation}\label{NS}
\begin{cases}
\partial_t v+(u \cdot \nabla) v +\nu(-\Delta)^{\beta}v +\nabla p =0 \\
v=u+(-\Delta)^{\alpha}u\\
\nabla \cdot u=\nabla \cdot v=0
\end{cases}
\end{equation}
and proved that this system is well posed when $\alpha+\beta \ge\frac 54$
(with $\alpha,\beta\ge0$); 
even some logarithmic corrections can be included, but we do not 
specify this detail, since it is not related to our analysis. It is worth mentioning the result of the current authors with Barbato
 in \cite{BBF}, where a stochastic version of the associated inviscid system to \eqref{NS} (when $\nu=0$) has been studied. In fact, by choosing an appropriate stochastic perturbation of the system to be formally conservative, they were able to prove global existence and uniqueness of solutions in law for $\alpha>\frac{3}{4}$. This is a very strong result although the uniqueness has to be understood in law. \\
Similar regularization have been used for the MHD models, see e. g.
\cite{Kazuo} and the references therein. Since these models are quite
different from the ones considered in the current paper, we don't
state their results and we refer interested readers to the literature
related to these models.
The physical motivation of these regularization defined in terms of
smoothing kernels is related to
a sub-grid length scale in the model and these kernels work as a kind
of filter with certain widths. An extensive explanation of these models can be
found in \cite{ot} and the references therein. 

Inspired by \cite{BMR},  in this paper we consider the modified 
Boussinesq system for $d=3$,
where the equation for the velocity has fractional dissipation 
whereas the temperature equation has no dissipation term; a
Leray-regularization for the velocity  appears in the quadratic terms.
This system is 
\begin{equation}\label{syst:iniz}
\begin{cases}
\partial_t v+(u \cdot \nabla) v +\nu(-\Delta)^{\beta}v +\nabla p =\theta e_3 \\
\partial_t \theta+u\cdot \nabla \theta=0\\
v=u+(-\Delta)^{\alpha}u\\
\nabla \cdot u=\nabla \cdot v=0
\end{cases}
\end{equation}
As in \cite{ot}, we work on a box and assume periodic boundary conditions.

Inspired by Ye \cite{Ye} and Barbato, 
Morandin, Romito \cite{BMR}, our goal is to prove 
 well posedness of system \eqref{syst:iniz} 
for $\alpha+\beta= \frac 54$ when $v_0, \theta_0$ are 
regular enough. So the interesting case is 
for $\beta<\frac 54$;  
indeed, the result of Ye corresponds 
to $\alpha=0$ and $\beta \ge \frac 54$ and that of Barbato, 
Morandin, Romito does not include the temperature equation, i.e. 
corresponds to our system \eqref{syst:iniz} with  $\theta_0=0$. 
We have to point out that  the temperature satisfies
a  pure transport equation, without thermal diffusivity; hence, the uniqueness 
result for the unknown $\theta$ requires the velocity to be smooth enough and this imposes
$\beta$ to be not too small.
We point out that \cite{KC} and \cite{Se} deal with some regularized 
Boussinesq system similar to \eqref{syst:iniz}; however they consider an easier case, since they set 
$\alpha=\beta=1$,  the bilinear term in the equation for the velocity is $(u\cdot \nabla) u$ instead of our $(u\cdot \nabla)v$
and the equation for the temperature is dissipative, i.e. there is the term $-\kappa \Delta \theta$ in the l.h.s..

We can summarize our result in the following
\begin{thm}
Assume $\frac 12<\beta<\frac 54$ with
\[
\alpha+\beta = \frac 54.
\] 
Then, system \eqref{syst:iniz} has a unique global smooth solution for any 
smooth initial conditions $v_0, \theta_0$.
\end{thm}

\begin{rmk}
Notice that assuming $\alpha+\beta=\frac 54$, the condition $\frac 12<\beta<\frac 54$
 is equivalent to $0<\alpha<\frac 34$. 
Let us notice that our techinque works also in the easier case
$\alpha+\beta >\frac 54$ with  $\beta>\frac 12$, $\alpha\ge 0$
(see Remark \ref{oss-su->}); but
the result of  \cite{Ye}  for $\alpha=0$ and $\beta=\frac 54$ cannot
be obtained with our techinque.
\end{rmk}

Our proofs rely on the commutator estimates introduced in \cite{kp}, also used in \cite{Ye}.
However  in contrast to \cite{Ye},  we first prove global existence 
(for any $\alpha\ge 0$ and $\beta>0$) and then uniqueness 
of these solutions; moreover we consider different order of 
space regularity for $v$ and $\theta$ 
($H^r$-regularity for $v$ and $H^{r-\beta}$-regularity 
for $\theta$), whereas in \cite{Ye}  the same order of regularity for both $v$ and $\theta$ is considered.
We point out that the requirement on the regularity on the initial data is 
needed  only to guarantee uniqueness.

The paper is organized as follows.
Section 2 is devoted to the mathematical framework.  Our main  functional spaces, the  regularization operator 
$\Lambda^s$ with its properties given in Lemma \ref{prod} are defined. 
The bilinear operator of the Navier-Stokes equations, the transport 
operator and the commutator operator are defined and their properties 
are stated in Lemma \ref{B}, Remark \ref{oss-su-def} and Remark 
\ref{oss-regol} and Lemma \ref{comm}.
The main system  is then written in its abstract (operator) form and the 
definition of weak solutions is given. 
At the end of this section, we recall the Gagliardo-Nirenberg and 
Br\'ezis-Gallouet-Wainger inequalities and  some continuity results.  
In Section 3, we prove global existence of weak solutions with their 
uniform estimates. Slightly better estimates are performed. However, 
they are not enough to prove the uniqueness of solutions. The main result 
of the paper is stated in Section 4, Theorem \ref{Tr}, where we prove global 
existence of regular solutions;  this regularity is enough to prove 
uniqueness of solutions and their continuous dependence with respect to the 
initial conditions, see Theorem \ref{uniq} and Theorem \ref{cc}.
 Let us point out that the results of Section 4 provide
 Theorem 1, i.e.
every smooth initial data gives rise to a 
unique smooth solution. Section 5 is devoted to showing in more details 
the crucial estimates  used  in Section 4.

\section{Mathematical framework}
We consider the evolution for positive times and the spatial variable belongs 
to a bounded domain  of $\mathbb R^3$; for simplicity and because
of the lack of natural boundary conditions,
we work on the torus, i.e. the spatial variable $x \in
\mathbb T=[0,2\pi]^3$  and  periodic boundary conditions are assumed. 
We set $L_p=L^p(\mathbb T)$.

As usual in the periodic setting, 
we can restrict ourselves to deal with initial data  with 
vanishing spatial averages; then the solutions will enjoy the same property
 at any fixed time $t>0$.

Therefore we can represent any $\mathbb T$-periodic 
function $f :\mathbb R^3 \to \mathbb R$
as
\[
f(x)=\sum_{k \in \mathbb Z^3_0}f_k e^{i k\cdot x}, \qquad
{\rm with } \;f_k \in \mathbb C, \ f_{-k}=\overline f_k\quad \forall k 
\]
where $\mathbb Z^3_0=\mathbb Z^3\setminus 0$.
For $s \in \mathbb R$ we define the spaces
\[
H^s=\{f=\sum_{k \in \mathbb Z^3_0}f_k e^{i k\cdot x}:
f_{-k}=\overline f_k \;{\rm and }\;
\sum_{k \in \mathbb Z^3_0} |f_k|^2|k|^{2s}<\infty\} .
\]
They are a Hilbert spaces with scalar product
\[
\langle f,g\rangle_{H^s}=\sum_{k \in \mathbb Z^3_0} f_k g_{-k}|k|^{2s} .
\]
We simply denote by $\langle f,g\rangle$ the  scalar product in $H^0$
and also the dual pairing of $H^s-H^{-s}$, i.e.
 $\langle f,g\rangle=\sum_k f_k g_{-k}$.

The space $H^{s+\epsilon}$ is compactly embedded in $H^s$ for any $\epsilon>0$.
Moreover, we recall the Sobolev embeddings:
if $0\le s< \frac{3}{2}$ and 
$\frac{1}{p}=\frac{1}{2}-\frac{s}{3}$, then 
$H^{s}\subset L_p$ and  there exists a constant $C$ (depending on $s$ and $p$)
such that
\begin{equation}
\|f\|_{L_p}\le C \|f\|_{H^{s}}.
\end{equation}
If $s=\frac 32$, then 
\[
 \|f\|_{L_p}\le C \|f\|_{H^s} \qquad \text{ for any finite } p
\]
and if $s>\frac 32$, then 
\[
 \|f\|_{L_\infty}\le C \|f\|_{H^s}.
\]
We shall often use the following inequality, merging the  two latter ones:
\begin{equation}\label{Sob-32}
\text{ if } s\ge \frac 32 \;\text{ then }
\|f\|_{L_p}\le C \|f\|_{H^s} \qquad \text{ for any finite } p
\end{equation}
Hereafter, we denote by the same symbol $C$ different constants.

Similarly, we define the spaces for the divergence free velocity vectors, which are 
periodic and have zero spatial average.
For $w:\mathbb R^3\to \mathbb R^3$ we write formally 
\[
w(x)=\sum_{k \in \mathbb Z^3_0}w_k e^{i k\cdot x}, \qquad
{\rm with }\; w_k \in \mathbb C^3, \ w_{-k}=\overline w_k, w_k \cdot k=0 \quad \forall k
\]
and for $s \in \mathbb R$  define
\[
V^s=\{w=\sum_{k \in \mathbb Z^3_0}w_k e^{i k\cdot x}:  w_{-k}=\overline w_k,\ 
 w_k \cdot k=0 \;{\rm and }\;
\sum_{k \in \mathbb Z^3_0} |w_k|^2|k|^{2s}<\infty\} .
\]
This is a Hilbert space with scalar product
\[
\langle v,w\rangle_{V^s}=\sum_{k \in \mathbb Z^3_0} v_k\cdot  w_{-k}|k|^{2s} .
\]

We define the linear operator $\Lambda=(-\Delta)^{1/2}$, i.e.
\[
f=\sum_{k \in \mathbb Z^3_0}f_k e^{i k\cdot x}
\Longrightarrow
\Lambda f=\sum_{k \in \mathbb Z^3_0}|k| f_k e^{i k\cdot x}
\]
and its powers $\Lambda^s$: 
$\Lambda^s f=\sum_{k \in \mathbb Z^3_0}|k|^s f_k e^{i k\cdot x}$;
hence $\Lambda^2=-\Delta$.
Note, in particular that $\Lambda^{s}$ maps 
$H^{r}$ onto $H^{r-s}$.
\\
For simplicity, we shall use the same notation for $\Lambda$ in the scalar 
spaces $H^s$ and in the vector spaces $V^s$. 

Let us denote by $\Pi$ the Leray-Helmholtz projection from $L_2$ 
onto $V^0$. The operators $\Pi$ and $\Lambda^s$ commute.

Finally we define the  bilinear operator $B: V^1\times V^1\to V^{-1}$ by
\[
\langle B(u,v),w\rangle=\int_{\mathbb T} \Big((u\cdot \nabla)v\Big)\cdot w \ dx
\]
i.e.
$ B(u,v)=\Pi\Big( (u\cdot \nabla)v\Big)$ for smooth vectors $u,v$.

We summarize the properties of the nonlinear terms; these are classical results, 
see e.g. \cite{temamP}.
\begin{lemma}\label{B}
For any $u,v,w\in V^1$ and $\theta,\eta \in H^1$ we have
\begin{equation}\label{tril0}
\langle B(u,v), w\rangle =-\langle B(u,w), v\rangle, \quad \langle B(u,v), v\rangle =0,
\end{equation}
\begin{equation}\label{tril1}
\langle v \cdot \nabla \theta, \eta\rangle =-\langle v \cdot \nabla \eta, \theta\rangle, \qquad \langle v \cdot \nabla \theta, \theta\rangle =0
\end{equation}
\eqref{tril0} holds more generally
for   any $u,v,w$ giving a meaning to the trilinear forms,
as stated precisely in the following:
\begin{equation}\label{stimeB}
\langle B(u,v),w\rangle \le C\|u\|_{V^{m_1}} \|v\|_{V^{1+m_2}} \|w\|_{V^{m_3}}
\end{equation}
with the non negative parameters fulfilling 
\[
m_1+m_2+m_3\ge \frac 32 \;\qquad\text{ if } m_i\neq \frac 32  \text{ for any } i
\]
or 
\[
m_1+m_2+m_3>\frac 32 \qquad\qquad \text{ if } \  \exists \ m_i=\frac 32.\qquad\qquad\;
\]
\end{lemma}

Now, we are ready to give the abstract formulation of problem 
\eqref{syst:iniz}; 
we apply the projection operator $\Pi$ to the first equation in order 
to get rid of the pressure. In addition due to the periodic setting, we regularize $u$ in a different, 
but equivalent way.
Therefore, our system in abstract form is 
\begin{equation}\label{abc}
\begin{cases}
\partial_t v+B(u,v)  +\nu\Lambda^{2\beta}v =\Pi (\theta e_3) \\
\partial_t \theta+u\cdot \nabla \theta=0\\
v=\Lambda^{2\alpha}u
\end{cases}
\end{equation}

We focus our analysis on the unknowns $v$ and $\theta$. 
The pressure $p$ will be recovered 
by taking the curl of the equation for 
the velocity in \eqref{syst:iniz}, i.e. $p$ solves the equation
$\Delta p =-\nabla\cdot[(u \cdot \nabla)v-\theta e_3]
=-\nabla\cdot[(( \Lambda^{-2\alpha}v)\cdot \nabla)v-\theta e_3]$. 
Moreover, 
the unknown $u$ is directly related to $v$. 

Therefore we give the following definition in terms of  $v$ and $\theta$ only. 
The finite time interval $[0,T]$  is fixed throughout the paper.

\begin{definition}\label{def-deb}
Let $\alpha\ge 0$ and $\beta>0$.
We are given $v_0\in V^0, \theta_0\in H^0$.
We say that the couple $(v,\theta)$ is a weak solution to system \eqref{abc} 
over the time interval $[0,T]$ if 
\begin{align*}
&v \in L^\infty(0,T;V^0)\cap L^2(0,T; V^\beta)\cap C_w([0,T];V^0)\\
&\theta \in L^\infty(0,T;H^0)\cap C_w([0,T];H^0)
\end{align*}
and, given any $\psi\in  V^{\frac 52}$, $\phi\in H^{\frac 52}$,  
they satisfy 
\begin{align}\label{var-v}
\langle v(t),\psi\rangle & -\int_{0}^{t}\langle B(u(s), \psi),v(s)\rangle ds
+\nu\int_{0}^{t}\langle \Lambda^{\beta}v(s), \Lambda^{\beta} \psi\rangle ds\nonumber\\
&=\langle v_0,\psi\rangle + \int_{0}^{t}\langle \theta(s) e_3,\psi \rangle ds
\end{align}
\begin{equation}\label{var-theta}
\langle \theta(t),\phi\rangle 
-\int_{0}^{t}\langle u(s)\cdot  \nabla \phi, \theta(s) \rangle ds 
=\langle \theta_0,\phi\rangle 
\end{equation}
for  every  $t\in [0,T]$.
\end{definition}

\begin{rmk}\label{oss-su-def}
In the weak formulations above,  the trilinear terms are well defined; indeed, 
if $0<2\alpha+\beta\le \frac 32$ 
\begin{equation}\label{stimaBab}
\begin{split}
|\langle B(u,\psi),v\rangle| 
&\le C
\|u\|_{V^{2\alpha}} \|\psi\|_{V^{\frac 52-2\alpha-\beta}} \|v\|_{V^{\beta}}  
    \qquad{\rm by }\;\eqref{stimeB} 
\\
&\le C \|v\|_{V^0} \|v\|_{V^\beta}\|\psi\|_{V^{\frac 52}} 
\end{split}\end{equation}
and  if $2\alpha+\beta> \frac 32$
\begin{equation}\label{stimaB-1}
\begin{split}
|\langle B(u,\psi),v\rangle| 
&\le C
\|u\|_{L_\infty} \|\nabla \psi\|_{L_2} \|v\|_{L_2}  
\qquad \text{ by H\"older inequality}
\\
&\le C \|u\|_{V^{2\alpha+\beta}} \|v\|_{V^0}\|\psi\|_{V^1} 
         \qquad{\rm since }\; V^{2\alpha+\beta}\subset L_\infty
\\
& 
\le  C \|v\|_{V^{\beta}} \|v\|_{V^0}\|\psi\|_{V^1} .
\end{split}\end{equation}
Similarly for the temperature: \\if  $0<2\alpha+\beta< \frac 32$ 
\begin{equation}\begin{split}
|\langle u\cdot  \nabla \phi, \theta \rangle|
&\le
\|u\|_{L_{p_1}}\|\nabla \phi\|_{L_{p_2}} \|\theta\|_{L_2} \\
&\le C \|u\|_{V^{2\alpha+\beta}} \|\nabla \phi\|_{H^{\frac 32-2\alpha-\beta}}
        \|\theta\|_{H^0} \\
&\le C \|v\|_{V^{\beta}} \|\phi\|_{H^{\frac 52}}\|\theta\|_{H^0}  
\end{split}
\end{equation}
where we used first the H\"older inequality with 
$\frac 1{p_1}=\frac 12-\frac {2\alpha+\beta} 3\in (0,\frac 12),
\frac 1{p_2}=\frac 12-\frac 1{p_1}$ and then 
the embedding theorems;\\
if  $2\alpha+\beta\ge \frac 32$ 
\begin{equation}\begin{split}
|\langle u\cdot  \nabla \phi, \theta \rangle|
&\le
\|u\|_{L_{4}}\|\nabla \phi\|_{L_{4}}\|\theta\|_{L_2}\\
& \le C  \|u\|_{V^{2\alpha+\beta}} \|\nabla \phi\|_{H^{\frac 34}}
        \|\theta\|_{H^0} \\
&\le C  \|v\|_{V^{\beta}} \|\phi\|_{H^{\frac 74}}\|\theta\|_{H^0}  \\
&\le C  \|v\|_{V^{\beta}} \|\theta\|_{H^0} \|\phi\|_{H^{\frac 52}} 
\end{split}
\end{equation}
where we used first the H\"older inequality and then 
the embedding theorems $V^{2\alpha+\beta}\subset L_{q}$ for any finite $q$,
$H^{\frac 52}\subset H^{\frac 74}$, $H^{\frac 34}\subset L_4$.

For more regular solutions, the trilinear term 
$\langle  B(u,\psi),v\rangle$ is equal to
$-\langle  B(u,v), \psi\rangle$ and we recover 
the term appearing in the equation for the velocity. The same holds
for the temperature.
\end{rmk}

\begin{rmk}\label{oss-regol}
We point out that the estimates by means of Sobolev embeddings need
some restriction for the parameters; but, for bigger values of the parameters
they are easier to prove and the details 
will be skipped. This means for instance that \eqref{stimaBab} 
with \eqref{tril0} gives 
\[
\|B(u,v)\|_{V^{-\frac 52}}\le C \|v\|_{V^0} \|v\|_{V^\beta}
\]
assuming $2\alpha+\beta\le \frac 32$, whereas for $2\alpha+\beta> \frac 32$ 
we get something stronger in \eqref{stimaB-1}:
\[
\|B(u,v)\|_{V^{-1}}\le C \|v\|_{V^0} \|v\|_{V^\beta}
\]
which is proven in another way. But for sure, from the proof 
of \eqref{stimaBab} 
 one can say that 
$\|B(u,v)\|_{V^{-\frac 52}}\le C \|v\|_{V^0} \|v\|_{V^\beta}$ also for
$2\alpha+\beta> \frac 32$ without proving it.
\end{rmk}

In this last part of the section, we summarize the 
technical tools to be used later on.

To estimate an $L_\infty$-norm we use either 
the embedding theorem $H^{r}\subset L_\infty$ 
with $r>\frac 32$ or
 the  {\bf Br\'ezis-Gallouet-Wainger inequality} (see \cite{BG,BW}):\\
for any $r>\frac 32$ there exists a constant $C$ such that
\begin{equation}
\|g\|_{L_\infty}\le C  \|g\|_{H^{\frac 32}} 
\left( 1+\sqrt{\ln (1+\frac{\|g\|_{H^r}}{\|g\|_{H^{\frac 32}}}})  \right) .
\end{equation}
Actually, we shall use the stronger form of this inequality, 
as given for instance in \cite{Ye}:
for any $r>\frac 32$ there exists a constant $C$ such that
\begin{equation}\label{stimaBG}
 \|g\|_{L_\infty}
  \le 
  C\left(1+\|g\|_{H^{\frac 32}}  +\|g\|_{H^{\frac 32}} \ln(e+\|g\|_{H^r})\right) .
\end{equation}

\noindent
{\bf Gagliardo-Niremberg inequality} (see \cite{Ni})\\
Let $1 \le q,r \le \infty$, $0< s< m $,  $\frac sm \le a <1$ and 
\[
\frac 1p=\frac s3+\left(\frac 1q -\frac m3\right)a +\frac {1-a}r
\]
then there exists a constant $C$ such that
\begin{equation}
\|\Lambda^s g \|_{L_p}\le C \|g\|_{L_r}^{1-a}\|\Lambda^m g\|^a_{L_q} .
\end{equation}

\medskip
Define the commutator 
\[
[\Lambda^s,f]g= \Lambda^s(fg) - f \ \Lambda^s g .
\]
From \cite{kp}, \cite{kpv} we have

\begin{lemma}[Commutator lemma]\label{comm}
Let $s>0$, $1<p<\infty$  and $p_2,p_3 \in (1,\infty)$ be such that
\[
\frac 1p\ge  \frac 1{p_1}+\frac 1{p_2},\qquad
\frac 1p\ge \frac 1{p_3}+\frac 1{p_4}.
\]
Then
\[
\|[\Lambda^s,f]g\|_{L_p}\le C\left(
\|\nabla f\|_{L_{p_1}}\|\Lambda^{s-1}g\|_{L_{p_2}}+
\|\Lambda^s f\|_{L_{p_3}}\|g\|_{L_{p_4}}
\right) .
\]
\end{lemma}
and 
\begin{lemma}\label{prod}
Let $s>0$, $1<p<\infty$  and $p_2,p_3 \in (1,\infty)$ be such that
\[
\frac 1p\ge \frac 1{p_1}+\frac 1{p_2},\qquad
\frac 1p\ge \frac 1{p_3}+\frac 1{p_4}.
\]
Then
\[
\|\Lambda^s(fg)\|_{L_p} \le C\left(
\|f\|_{L_{p_1}}\|\Lambda^s g\|_{L_{p_2}}+
\|\Lambda^s f\|_{L_{p_3}}\|g\|_{L_{p_4}}
\right) .
\]
\end{lemma}

We shall use the commutator acting also on vectors; in particular
for $u,v \in \mathbb R^3, \theta \in \mathbb R$
\[
[\Lambda^s,u]\cdot \nabla \theta
=
\Lambda^s(u \cdot \nabla \theta)
-u \cdot \nabla \Lambda^s \theta
\]
and
\[
[\Lambda^s,u]\cdot \nabla v
=
\Lambda^s\big((u \cdot \nabla)v\big)
-(u \cdot \nabla) \Lambda^s v .
\]
Therefore
\begin{equation}\label{tril-con-comm1}
\langle \Lambda^s(u \cdot \nabla \theta), \Lambda^s \theta\rangle =
\langle [\Lambda^s,u]   \cdot \nabla \theta, \Lambda^s \theta\rangle +
\underbrace{\langle u \cdot \nabla \Lambda^s \theta , \Lambda^s \theta\rangle }
_{=0\text { by } \eqref{tril1}}
\end{equation}
and
\begin{equation}\label{tril-con-comm2}
\langle \Lambda^s\big((u \cdot \nabla)v \big), \Lambda^s v\rangle =
\langle [\Lambda^s,u]  \cdot \nabla v, \Lambda^s v\rangle +
\underbrace{\langle (u \cdot \nabla)\Lambda^s v , \Lambda^s v\rangle }
_{=0\text { by } \eqref{tril0}}
\end{equation}

About the continuity in time, we have the strong continuity result
(see \cite{St} or Lemma 1.4, Chap III in \cite{temam})
\begin{lemma}\label{lemmaC0} Let $s \in \mathbb R$ and $h>0$.\\
If $v \in L^2(0,T;V^{s+h})$ and $\frac {dv}{dt} \in L^2(0,T;V^{s-h})$, then
$v \in C([0,T];V^s)$ and
\[
\frac {d}{dt}\|v(t)\|_{V^s}^2
=
2 \langle \Lambda^{-h}\frac {dv}{dt}(t), \Lambda^{h}v(t)\rangle_{V_s} 
\]
\end{lemma}
and the weak continuity result (see \cite{St}).

\begin{lemma}\label{lemmaCw}
Let $X$ and $Y$ be Banach spaces, $X$ reflexive,
$X$ a dense subset of $Y$ and the inclusion map of $X$ into $Y$ continuous.
Then 
\[
L^\infty(0,T;X)\cap C_w([0,T];Y)= C_w([0,T];X).
\]
\end{lemma}

\section{Existence of weak solutions}
Existence of a global weak solution of system \eqref{abc} can be obtained easily; 
the technique is very similar to that for the classical Boussinesq system. 
The equation for $\theta$ is a pure transport equation; then
the $L_q$-norm of $\theta$ is conserved in time (for any $q\le +\infty$).
On the other hand, it is enough to have some regularization 
in the velocity equation (i.e. $\beta>0$) 
in order to get a weak solution as in Definition \ref{def-deb}; moreover, 
this solution satisfies an energy inequality.
Of course, the bigger are the parameters $\alpha, \beta$, the more regular is the velocity $v$.
\begin{thm}\label{teo0}
Let $\alpha\ge 0, \beta>0$ and $2\le q \le \infty$.
For any  $v_0\in V^0, \theta_0 \in L_q$, there exists a 
weak solution $(v, \theta)$ of \eqref{abc} on the time interval $[0,T]$.
Moreover
$$
 \theta \in C_{w}(0,T;L_q).
$$
\end{thm}
\pf
We define the finite dimensional projector operator $\Pi_n$ in $V^0$ as
$\Pi_n v=\sum_{0<|k|\le n} v_k e^{i k \cdot x}$ 
for $v=\sum_{k \in \mathbb Z^3_0} v_k e^{i k \cdot x}$; similarly for the scalar case,
 i.e. $\Pi_n$ in $H^0$. We set $B_n(u,v)=\Pi_n B(u,v)$.

We consider the finite dimensional approximation of system \eqref{abc}  
in the unknowns
$v_{n}=\Pi_n v$,  $u_{n}=\Pi_n u$ and $\theta_{n}=\Pi_n\theta$. 
 This is the Galerkin approximation for $n=1,2, \dots$
\begin{equation}
\begin{cases}
\partial_t v_n+B_n(u_n,v_n)  +\nu\Lambda^{2\beta}v_n =\Pi (\theta_n e_3) \\
\partial_t \theta_n+\Pi_n(u_n\cdot \nabla \theta_n)=0\\
v_n=\Lambda^{2\alpha}u_n
\end{cases}
\end{equation}

We take the $L_2$-scalar product of  the equation for the velocity $v_n$ with  $v_n$ itself; 
bearing in mind \eqref{tril0} we get
\[
\begin{split}
\frac 12 \frac{d}{dt}\|v_n(t)\|_{V^0}^2 +\nu\|v_n(t)\|^2_{V^\beta}
&=-\langle B_n(u_n(t), v_n(t)),v_n(t)\rangle+\langle \Pi(\theta_n(t) e_3) , v_n(t)\rangle  \\
&=-\langle B(u_n(t), v_n(t)),v_n(t)\rangle+\langle \theta_n(t) e_3 , v_n(t)\rangle  \\
&\le \frac 12 \|\theta_n(t)\|_{H^0}^2 +\frac 12 \|v_n(t)\|_{V^0}^2
\end{split}
\]
and similarly for the second equation
\[\begin{split}
\frac{d}{dt}\|\theta_n(t)\|_{H^0}^2
&=-\langle \Pi_n (u_n(t) \cdot \nabla \theta_n(t)),\theta_n(t)\rangle\\
&= -\langle u_n(t) \cdot \nabla \theta_n(t),\theta_n(t)\rangle=0 .
\end{split}\]
In both cases the trilinear forms vanish according to 
\eqref{tril0}, \eqref{tril1}.

Adding these estimates, by means of Gronwall's lemma we get the basic 
$L_2$-energy estimate: there exists a constant $K_1$ independent of $n$ such that
\[
\sup_{0\le t\le T}(\|v_n(t)\|_{V^0}^2+\|\theta_n(t)\|_{H^0}^2)
+\nu \int_0^T\|v_n(t)\|^2_{V^\beta}dt\le K_1
\]
for any $n$.

From the equation for the velocity $v_n$, one has that $\frac{dv_n}{dt}$ is expressed 
as the sum of three terms involving $v_n$, $u_n$ 
and $\theta_n$. In particular, the dissipative term 
$\Lambda^{2\beta}v_n\in L^2(0,T;V^{-\beta})$; by \eqref{stimaBab},
\eqref{stimaB-1}
we have $B_n(u_n,v_n)\in L^2(0,T;V^{-s})$ for some finite $s\ge 1$.
Therefore there exist constants $\gamma>0$ and $K_2$
independent of $n$,  such that
\begin{align*}
\|\frac{dv_n}{dt}\|_{L^2(0,T;V^{-\gamma})}^2\le K_2.
\end{align*}

This means that $v_{n}$ is bounded in 
$L^2(0,T; V^{\beta})\cap W^{1,2}(0,T; V^{-\gamma})$ 
(with $\beta>0$ and $\gamma>0$),
 which is compactly embedded in $L^2(0,T; V^{0})$ 
(see Lemma 2.2. in \cite{temam}). 
Hence we can extract a subsequence, still denoted by $\{v_{n}\}$ and 
$\{\theta_{n}\}$, such that
$$
v_{n}\longrightarrow v \quad {\rm weakly \ in }\quad L^2(0,T; V^{\beta})
$$
$$
v_{n}\longrightarrow v \quad {\rm weakly}^*\  {\rm in } \quad L^{\infty} (0,T; V^{0})
$$
$$
v_{n}\longrightarrow v \quad {\rm strongly\ in }\quad L^2(0,T; V^{0})
$$
$$
\theta_{n}\longrightarrow \theta \quad {\rm weakly}^* \  {\rm in }\quad 
    L^{\infty} (0,T; H^{0}).
$$
  
Using these convergences, it is a classical result
 to pass to the limit in the variational formulation \eqref{var-v} and 
\eqref{var-theta} and prove that $(v,\theta)$ is solution of \eqref{abc} and inherits  all the regularity from 
$(v_{n},  \theta_{n})$, i.e. 
\[
v \in L^\infty(0,T;V^0)\cap L^2(0,T; V^\beta),\qquad
\theta \in L^\infty(0,T;H^0).
\]

Moreover, it is a classical result (see \cite{Ye}) that 
\begin{equation}\label{stima-base-theta_n}
\sup_{0\le t\le T}\|\theta_{n}(t)\|_{L_q}\le \|\theta_0\|_{L_q}
\end{equation}
for any $q \le \infty$. 

Hence, the sequence $\left\{\theta_{n}\right\}_n$ is uniformly bounded in 
$L^{\infty} (0,T; L_q)$ which implies (up to a subsequence  still denoted $\theta_{n}$) that
 \[
\theta_{n}\longrightarrow \theta \quad {\rm weakly}^*\ {\rm in }
\quad L^{\infty} (0,T; L_{q})
\]
and 
\begin{equation}\label{stima-base-theta}
\sup_{0\le t\le T}\|\theta(t)\|_{L_q}\le \|\theta_0\|_{L_q}.
\end{equation}

Now, let us prove that $v \in C_w([0,T];V^0)$ and 
$\theta \in C_{w}([0,T];L_q)$. 
\\
We integrate in time the equation for $v$:
\[
v(t)=v_0+\int_0^t [-B(u(s),v(s))-\nu \Lambda^{2\beta}v(s)+\Pi\theta(s)e_3] ds .
\]
Bearing in mind \eqref{tril0} and 
the estimates of Remark \ref{oss-su-def}, we get that
$B(u,v)\in L^2(0,T;V^{-\frac 52})$; therefore $v \in C([0,T];V^{-m})$ 
for some positive $m$.
By Lemma \ref{lemmaCw} we get that $v \in C_w([0,T];V^0)$.

Now we look for  the weak continuity of $\theta$.
Assume that $\phi\in C^\infty_{\#}(\mathbb T)$ 
which is the space of $C^\infty$ functions on $\mathbb T$ that are periodic. 
Then for $t,s\in [0,T]$, we have that
\begin{align*}
|\langle \theta(t)-\theta(s),\phi\rangle|
&=|\int_{s}^{t}\langle u(r)\cdot\nabla \phi, \theta(r) \rangle dr|
\\
&\le \int_{s}^{t}\|\nabla\phi\|_{L_\infty} \|u(r)\|_{L_2} \|\theta(r)\|_{L_2}dr\\
&\le \|\nabla\phi\|_{L_\infty}  \|\theta\|_{L^\infty(0,T; H^0)} 
\int_s^t \|u(r)\|_{V^0} dr.
\end{align*}
Using the density of $C^\infty_{\#}(\mathbb T)$ in $L_{q^\prime}$ (with $\frac 1q +\frac 1{q^\prime}\le 1$), we deduce that
\[
\lim_{t \to s} \langle \theta(t)-\theta(s),\phi\rangle = 0\quad \forall \phi\in L_{q^\prime}
\]
which means that $\theta \in C_{w}([0,T];L_q)$.  A similar argument can be used for $q=\infty$ and this completes the proof.
\hfill $\Box$

\begin{rmk}
Take $\alpha\ge 0$ and $\beta>0$ such that
\[
2\alpha+\beta\le \frac 32 , \qquad
\alpha+\beta\ge \frac 54.
\]
For this to hold it is necessary that  $\alpha $ is not too big ($\alpha\le \frac 14$) 
and $\beta$ not too small ($1\le \beta\le \frac 32$).
Then,  from the first estimate in 
\eqref{stimaBab} we get $B(u,v)\in L^2(0,T;V^{-\beta})$. 
Hence, going back to the proof of the previous theorem we get
that $\frac{dv}{dt}\in L^2(0,T;V^{-\beta})$; by Lemma \ref{lemmaC0} this implies 
that $v \in C([0,T];V^0)$, which is stronger than the weak continuity result of Theorem
\ref{teo0} (see Definition 3).
\end{rmk}

In addition, for more regular initial data we have

\begin{thm}[More regularity]\label{Ts=beta}
We are given parameters $\alpha$ and $\beta$  with $\frac 12 <\beta<\frac 54$ and 
\[
\alpha+\beta= \frac 54.
\]
Then, given  $v_0\in V^\beta, \theta_0 \in H^0$, 
any weak  solution of \eqref{abc} obtained in  Theorem 
\ref{teo0} is more regular; indeed, the velocity is more regular
\[
v \in C([0,T];V^\beta) \cap L^2(0,T;V^{2\beta}).
\]
\end{thm}
\pf
We look for a priori estimates for $v$.
We proceed as before, but for more regular norms. We have  
\[
\begin{split}
\frac 12 \frac{d}{dt}\|v(t)\|_{V^\beta}^2 &+\nu\|v(t)\|^2_{V^{2\beta}}
\\
&=-\langle B(u(t), v(t)),\Lambda^{2\beta}v(t) \rangle
  +\langle \Pi \theta(t) e_3 , \Lambda^{2\beta} v(t)\rangle  \\
& \le |\langle\Lambda^\beta\big((u (t)\cdot\nabla) v(t)\big),\Lambda^{\beta}v(t)\rangle|
  +\|\theta(t)\|_{L_2} \| \Lambda^{2\beta} v(t)\|_{L_2}\\
& \le \|[\Lambda^\beta,u(t) ] \cdot\nabla v(t)\|_{L_2} \|\Lambda^{\beta}v(t)\|_{L_2}
  + \|\theta(t)\|_{H^0} \|v(t)\|_{V^{2\beta}}
\end{split}
\]
where we used \eqref{tril-con-comm2}.

We use the Commutator Lemma \ref{comm}
\begin{equation}\label{stime-Comm34}\begin{split}
\|[\Lambda^\beta,u ] \cdot\nabla v\|_{L_2} 
&\le C \left( \|\Lambda u\|_{L_{P_1}} \| \Lambda^{\beta}v\|_{L_{p_2}}+
\|\Lambda^{\beta} u\|_{L_{p_3}}\|\nabla v\|_{L_{p_4}}\right) 
\\
&=
C\left( \|\Lambda^{1-2\alpha} v\|_{L_{p_1}} \|\Lambda^{\beta} v\|_{L_{p_2}}+
\|\Lambda^{\beta-2\alpha} v\|_{L_{p_3}}\|\nabla v\|_{L_{p_4}}\right) 
\end{split}
\end{equation}
and we want to estimate further with $C\|v\|_{V^\beta}\|v\|_{V^{2\beta}}$.

For this we take 
\[
\frac 1{p_1}=\frac 12-\frac{\beta-1+2\alpha}3 \equiv \frac \beta3,
 \qquad \qquad
\frac 1{p_2}=\frac 12-\frac \beta 3 .
\]
The assumption $\frac 12 <\beta<\frac 54$ provides $3<p_2<12$
and by Sobolev embedding
\[
\|\Lambda^{1-2\alpha} v\|_{L_{p_1}}\le C \|v\|_{V^\beta},\qquad 
\|\Lambda^\beta  v\|_{L_{p_2}}\le C \|v\|_{V^{2\beta}} .
\]

For the latter two terms in \eqref{stime-Comm34} we choose 
\[
\frac 1{p_3}=\frac 12-\frac {2\alpha}3,
\qquad
\frac 1{p_4}=\frac 12-\frac 1{p_3}\equiv\frac 12 -\frac{2\beta-1}3 .
\]
The assumption $0 <\alpha<\frac 34$, i.e.
$\frac 12 <\beta<\frac 54$, 
provides $2<p_3<\infty$ and by Sobolev embedding
\[
\|\Lambda^{\beta-2\alpha} v\|_{L_{p_3}}\le C \|v\|_{V^{\beta}},
\qquad
\|\nabla v\|_{L_{p_4}}\le C\|v\|_{V^{2\beta}} .
\]

Hence, we conclude that
\begin{equation}\label{st-vH}
\begin{split}
\frac 12 \frac{d}{dt}&\|v(t)\|_{V^\beta}^2 +\nu\|v(t)\|^2_{V^{2\beta}}\\
&\le
C \|v(t)\|_{V^\beta} \|v(t)\|_{V^{2\beta}} \|v(t)\|_{V^{\beta}} +
C \|\theta(t)\|_{H^0} \|v(t)\|_{V^{2\beta}}
\\& \le
\frac \nu 2 \|v(t)\|^2_{V^{2\beta}}+
C_\nu \| v(t)\|_{V^\beta}^4  + C_\nu \| \theta(t)\|^2_{H^0}
\end{split}
\end{equation}
by Young inequality. In particular,
\[
\frac{d}{dt}\|v(t)\|_{V^\beta}^2\le C_\nu \| v\|_{V^\beta}^4
+ C_\nu \| \theta\|^2_{H^0}.
\]
Since $v \in L^2(0,T;V^\beta)$ and $\theta \in L^\infty(0,T;H^0)$ 
from the previous theorem, 
we can proceed by means of Gronwall lemma to get the estimate for 
the $L^\infty(0,T;V^\beta)$-norm:
\[
\sup_{0\le t\le T} \|v(t)\|_{V^\beta}^2\le
\|v_0\|_{V^\beta}^2 e^{C_\nu\int_0^T\|v(s)\|^2_{V^\beta}ds}
+C_\nu \int_0^T e^{C_\nu\int_r^T\|v(s)\|^2_{V^\beta}ds}\| \theta(r)\|^2_{H^0}dr .
\]
Integrating in time \eqref{st-vH}, we also get 
\[%\begin{multline}
\frac \nu 2 \int_0^T \|v(t)\|^2_{V^{2\beta}}\ dt\le 
\frac 12\|v_0\|_{V^\beta}^2%\\
+C_\nu \| v\|_{L^\infty(0,T;V^{\beta})}^4
+ C_\nu \int_0^T \| \theta(t)\|^2_{H^0} dt .
\]%\end{multline}
Summing up, we get that 
$v \in L^\infty(0,T;V^\beta)\cap L^2(0,T;V^{2\beta})$.

Now, we study the time regularity.
We recall   property \eqref{stimeB} for  the nonlinear term $B(u,v)$
with 
$ m_1= 2\alpha>0$, $m_2=2\beta-1> 0$, $m_3=0$ (we are in the first case, with all $m_i\neq \frac 32$ and thus we take $m_1+m_2+m_3=\frac 32$). We have
\begin{align*}
\|\frac{dv}{dt} (t)\|_{L_2}
&\le  \| B(u(t),v(t))\|_{L_2}+\nu \|\Lambda^{2\beta}v(t)\|_{L_2}+
 \|\theta (t)e_{3}\|_{L_2}\\
&\le 
 C\|u(t)\|_{V^{2\alpha}}  \|v(t)\|_{V^{2\beta}}  
 +\nu\|v(t)\|_{V^{2\beta}} +\|\theta(t)\|_{H^0}\\
&=
  C \|v(t)\|_{V^0}\|v(t)\|_{V^{2\beta}}+
 \nu\|v(t)\|_{V^{2\beta}} +\|\theta(t)\|_{H^0}
\end{align*}
Hence, using the regularity of  $v,\theta$ we get that
\[
\frac{dv}{dt} \in L^{2}(0,T; V^{0}).
\]
Now using Lemma \ref{lemmaC0}, we deduce that $v\in C([0,T]; V^\beta)$.
 \hfill $\Box$

\begin{rmk}\label{oss-su->}
The result of Theorem \ref{Ts=beta} still holds true under the assumption that
$\alpha+\beta>\frac 54$ with $\beta>\frac1 2$. This is trivial when 
we add the condition  $\alpha>0$, since the framework is similar to (but
easier than) that in the above proof. 
So, it remains to consider the case $\alpha=0$ and $\beta>\frac 54$.
To estimate the r.h.s. in \eqref{stime-Comm34} we choose 
$p_1=\frac{12}5, p_2=12$, $p_3=2$, $p_4=\infty$ so to get
\[
\|\Lambda v\|_{L_{p_1}}\le C\|v\|_{V^{\frac 54}}\le C\|v\|_{V^\beta}
\]
\[
\|\Lambda^\beta v\|_{L_{p_2}} \le C \|v\|_{V^{\frac 54+\beta}}\le 
\|v\|_{V^{2\beta}}
\]
\[
\|\nabla v\|_{L_{p_4}} \le C\|v\|_{V^{2\beta}}
\]
In the study of the time regularity, we choose $m_1=m_3=0$ and
$m_2=2\beta-1>\frac 32$ and conclude as above.

Similar remarks hold for the proofs of the Appendix, which are still valid
when assuming $\alpha+\beta>\frac 54$ with $\beta>\frac1 2$.

However,
our technique requires $\beta>\frac 12$. This might be improved as in \cite{BMR}; this is postponed to future work.
\end{rmk}

%%%%%%%%%%%%%%%%%%%%%%%%%%%%%%%%%%%%%%%%%%%%%%%%%%%%%%%
\section{Regular solutions: global existence, uniqueness and continuous dependence on the initial data}
The regularity of solutions from the previous section is not enough to prove 
 uniqueness. 
To this end, we seek  classical  
solutions.  These are solutions for which the spatial derivatives
in the equations of \eqref{abc} exist. Indeed, we shall get that 
$v \in C([0,T];V^r)\cap L^2(0,T;V^{r+\beta})$ and 
$\theta \in  C([0,T];H^{r-\beta})$
with $r>\beta+ 1$. The crucial point is to show that these regular 
solutions are defined on any given time interval $[0,T]$; 
their local existence is easy to prove.

Unlike the previous section, here we will consider $H^s$-regularity for $\theta(t)$ (with $s>0$).
This will help prove the uniqueness of solutions.

\begin{thm}\label{Tr}
We are given non negative parameters with $\frac 12 <\beta<\frac 54$ and 
\begin{equation}
\alpha+\beta= \frac 54 . \label{cond-ab}
\end{equation}
Let 
\[
r>\max\left( 2\beta,\beta+1\right).
\] 
Then, 
for any $v_0\in V^r, \theta_0 \in H^{r-\beta}$, there exists a 
solution $(v,\theta)$ to \eqref{abc} such that 
\[
v \in C([0,T];V^r)\cap L^2(0,T;V^{r+\beta}), 
\qquad \theta \in C([0,T];H^{r-\beta}).
\]
\end{thm}
\pf
We proceed as before. 
We take the $L_2$-scalar product of the first equation of \eqref{abc}
 with $\Lambda^{2r}v$;
then
\begin{equation}\label{stima1}\begin{split}
\frac 12 \frac{d}{dt}&\|v(t)\|_{V^r}^2 +\nu\|v(t)\|^2_{V^{r+\beta}}
\\
&=-\langle B(u(t), v(t)),\Lambda^{2r}v(t)\rangle
  +\langle \theta(t) e_3 , \Lambda^{2r} v(t)\rangle  \\
&=-\langle B(\Lambda^{-2\alpha} v(t), v(t)),\Lambda^{2r}v(t)\rangle
+\langle \Lambda^{r-\beta}\theta(t) e_3 , \Lambda^{r+\beta} v(t)\rangle
\\
&\le 
 C \|v(t)\|_{V^{2\beta}} \|v(t)\|_{V^{r+\beta}} \|v(t)\|_{V^{r}} 
+C \|\theta(t)\|_{H^{r-\beta}} \|v(t)\|_{V^{r+\beta}}\\
&\le
\frac \nu 4 \|v(t)\|_{V^{r+\beta}}^2+C_\nu\|v(t)\|_{V^{2\beta}}^2\|v(t)\|_{V^{r}}^2+C_\nu
\|\theta(t)\|_{H^{r-\beta}}^2
\end{split}
\end{equation}
where we used first Lemma \ref{lv1} and then
Young inequality.

Now for  $\theta$, we take the $L_2$-scalar product of the second  equation of \eqref{abc}
with $\Lambda^{2r-2\beta}\theta(t)$; then
\[
\frac 12 \frac{d}{dt}\|\theta(t)\|_{H^{r-\beta}}^2 
=-\langle u(t) \cdot \nabla
\theta(t),\Lambda^{2r-2\beta}\theta(t)\rangle .
\]
We estimate the r.h.s.
\[
\begin{split}
\langle u \cdot \nabla \theta,\Lambda^{2r-2\beta}\theta\rangle
&=\langle\Lambda^{r-\beta}( u \cdot \nabla\theta), \Lambda^{r-\beta}\theta\rangle\\
&=\langle [\Lambda^{r-\beta}, u]\cdot\nabla\theta,\Lambda^{r-\beta}\theta\rangle \;\text{ by } \eqref{tril-con-comm1}\\
& \le \|[\Lambda^{r-\beta}, \Lambda^{-2\alpha}v] \cdot \nabla\theta\|_{L_2}
   \|\Lambda^{r-\beta}\theta\|_{L_2}\\
\intertext{and the Commutator Lemma \ref{comm} gives }
&\le C\left( \|\Lambda^{1-2\alpha} v\|_{L_\infty}\|\Lambda^{r-\beta}\theta\|_{L_2} + \|\Lambda^{r-\beta-2\alpha}v\|_{L_{q_3}}\|\Lambda \theta\|_{L_{q_4}}
 \right)\|\Lambda^{r-\beta}\theta\|_{L_2}\\
\intertext{with $\frac 1{q_3}+\frac 1{q_4}\le \frac 12$; 
we continue by means of the Br\'ezis-Gallouet-Wainger
 estimate \eqref{stimaBG} (with $g=\Lambda^{1-2\alpha} v$) and Lemma \ref{ltt}}
&
\le C\left(1+\|\Lambda^{\frac 52-2\alpha} v\|_{L_2}+\|\Lambda^{\frac 52-2\alpha} v\|_{L_2}  
\ln(e+\|v\|_{V^{m+1-2\alpha}})\right) \|\theta\|_{H^{r-\beta}}^2\\
& \qquad +C \|v\|_{V^{2\beta}}^{a} \|v\|_{V^{r+\beta}}^{1-a} 
\|\theta\|_{L_q}^{1-a} \|\theta\|_{H^{r-\beta}}^{1+a}\\
\intertext{for any $m>\frac 32$ 
and for suitable $q>2$, $a\in (0,1)$; $m$ will be chosen later on. 
Finally we use that $V^{2\beta}= V^{\frac 52-2\alpha}$ and $V^{m+1-2\alpha}=V^{m+2\beta-\frac 32}$,
since $\alpha+\beta= \frac 54$:}
&
\le C\left(1+\|v\|_{V^{2\beta}}+\| v\|_{V^{2\beta}}  
\ln(e+\|v\|_{V^{m+2\beta-\frac 32}})\right) \|\theta\|_{H^{r-\beta}}^2\\
& \qquad +C \|v\|_{V^{2\beta}}^{a} \|v\|_{V^{r+\beta}}^{1-a} 
\|\theta\|_{L_q}^{1-a} \|\theta\|_{H^{r-\beta}}^{1+a}.
\end{split}
\]

Now, we use 
Young inequality: 
\[
\|v\|_{V^{2\beta}}^{a} \|v\|_{V^{r+\beta}}^{1-a} \|\theta\|_{L_q}^{1-a} 
\|\theta\|_{H^{r-\beta}}^{1+a}
\le
\frac \nu 4 \|v\|_{V^{r+\beta}}^2+C_\nu \|v\|_{V^{2\beta}}^{\frac{2a}{1+a}}
\|\theta\|_{L_q}^{\frac{2(1-a)}{1+a}}\|\theta\|_{H^{r-\beta}}^{2} .
\]

Set  $\phi:=\|v\|_{V^{2\beta}}^{\frac{2a}{1+a}}\|\theta\|_{L_q}^{\frac{2(1-a)}{1+a}}$; then 
$\phi \in L^1(0,T)$ according to Theorem \ref{Ts=beta}  and \eqref{stima-base-theta}.
Thus
\begin{multline}\label{stima2}
\frac 12 \frac{d}{dt}\|\theta(t)\|_{H^{r-\beta}}^2\le
C\big(1+\|v(t)\|_{V^{2\beta}}+\|v(t)\|_{V^{2\beta}}  
\ln(e+\|v(t)\|_{V^{m+2\beta-\frac 32}})\big) \|\theta(t)\|_{H^{r-\beta}}^2\\+\frac \nu 4
 \|v(t)\|_{V^{r+\beta}}^2+C_\nu \phi(t) \|\theta(t)\|_{H^{r-\beta}}^{2} .
\end{multline}

Adding the estimates \eqref{stima1} for $v$ and \eqref{stima2} for 
$\theta$, we get
\begin{multline}
\frac{d}{dt}(\|v(t)\|_{V^r}^2+\|\theta(t)\|_{H^{r-\beta}}^2)
+\nu\|v(t)\|^2_{V^{r+\beta}}
\le
C \| v(t)\|_{V^{2\beta}}^2\|v(t)\|_{V^r}^2\\
+C \big(1+ \|v(t)\|_{V^{2\beta}}+\|v(t)\|_{V^{2\beta}} 
\ln(e+\|v\|_{V^{m+2\beta-\frac 32}})+\phi(t)\big)
\|\theta(t)\|_{H^{r-\beta}}^2.
\end{multline}
Recall that $r>2\beta$ by assumption; 
then there exists $m>\frac 32$ such that
$V^r\subset V^{m+2\beta-\frac 32}$.
Thus, we get
\begin{multline}\label{stv2r}
\frac{d}{dt}(\|v(t)\|_{V^r}^2+\|\theta(t)\|_{H^{r-\beta}}^2)
+\nu\|v(t)\|^2_{V^{r+\beta}}
\le
C  \| v(t)\|_{V^{2\beta}}^2\|v(t)\|_{V^r}^2\\
+C \Big(1+\|v(t)\|_{V^{2\beta}}+\|v(t)\|_{V^{2\beta}} \ln(e+\|v(t)\|_{V^{r}})+\phi(t)\Big)
\|\theta(t)\|_{H^{r-\beta}}^2 
\end{multline}
Set $X(t)=\|v(t)\|_{V^r}^2+\|\theta(t)\|_{H^{r-\beta}}^2$. Then, 
from \eqref{stv2r} we easily get
\[\begin{split}
\frac {dX}{dt} (t)&\le C\Big(1+\|v(t)\|_{V^{2\beta}} 
\ln(e+1+X(t))+\|v(t)\|^2_{V^{2\beta}}+\phi(t)\Big)X(t)
\\
&
\le C\Big(1+\|v(t)\|_{V^{2\beta}} \ln(e+1+X(t))
    +\|v(t)\|^2_{V^{2\beta}}+\phi(t)\Big)(e+1+X(t)).
\end{split}
\]
This implies that $Y(t)=\ln (e+1+X(t))$ satisfies
\[
Y^\prime(t)
\le C \Big(1+\|v(t)\|_{V^{2\beta}}  Y(t)+\|v(t)\|^2_{V^{2\beta}}+\phi(t)\Big).
\]
By Gronwall lemma we get 
$$
\sup_{0\le t\le T}Y(t) \le Y(0)e^{C\int_0^T \|v(s)\|_{V^{2\beta}}ds}
+C\int_0^T e^{C \int_s^T \|v(r)\|_{V^{2\beta}}dr} 
  \big(1+\|v(s)\|^2_{V^{2\beta}}+\phi(s)\big)ds .
$$
Since $v \in  L^2(0,T;V^{2\beta})$ by Theorem \ref{Ts=beta} and $\phi \in L^1(0,T)$, we get that 
\[
\sup_{0\le t\le T}Y(t) \le K_3
\]
and therefore going back to the unknown $X$ 
\[
\sup_{0\le t\le T}X(t) \le K_4;
\]
from \eqref{stv2r}, after integration on $[0,T]$ we get also
\[
\int_0^T \|v(t)\|_{V^{r+\beta}}^2 dt \le K_5.
\]

Therefore we have proved that
\[
v \in L^\infty(0,T;V^r)\cap L^2(0,T;V^{r+\beta}), 
\qquad \theta \in L^\infty(0,T;H^{r-\beta}).
\]
Now we consider the continuity in time.
Lemma \ref{prod} (with $p=p_2=2$, $p_1=\infty$)
gives
 \[
\|B(u,v)\|_{V^{r-\beta}}
\le C\left(\|u\|_{L_\infty}\|v\|_{V^{r-\beta+1}}
    +\|\Lambda^{r-\beta}u\|_{L_{p_3}} \|\Lambda v\|_{L_{p_4}}\right).
\]
By Sobolev embeddings we get
\[
\|u\|_{L_\infty}\le C \|\Lambda^{-2\alpha} v\|_{L_\infty}
\le C\|v\|_{V^r}
\]
since $r+2\alpha>\frac32$ (this comes from the assumption  $r>2\beta=\frac 52-2\alpha$),
and
\[
\|v\|_{V^{r-\beta+1}}\le C \|v\|_{V^{r+\beta}}
\]
since $\beta> \frac 12$.

Now we choose $p_3 \in (2,\infty)$ and $p_4$ such that $\frac 1{p_3}+\frac 1{p_4}=\frac 12$.
When $\beta>1$ we set $\frac 1{p_3}=\frac
12-\frac{\beta+2\alpha}3\equiv \frac{\beta-1}3$
and $\frac 1{p_4}=\frac 12-\frac {\beta-1}3$, 
so to get by Sobolev embedding
\begin{equation}\label{eq-11}
\|\Lambda^{r-\beta}u\|_{L_{p_3}}= \|\Lambda^{r-\beta-2\alpha} v\|_{L_{p_3}}
\le C \|v\|_{V^r}
\end{equation}
\[
\|\Lambda v\|_{L_{p_4}} \le C  \|v\|_{V^\beta}\le C  \|v\|_{V^{r+\beta}}\qquad \text{ for any } r\ge 0
\]
whereas when $\beta\le 1$ 
we have that 
\begin{equation}\label{nuovo-argomento}
\|\Lambda v\|_{L_{p_4}} \le C\|\Lambda v\|_{V^{r+\beta-1}}=C\|v\|_{V^{r+\beta}}
\end{equation}
for some $p_4\in (2,\infty)$ as soon as $r+\beta-1>0$
(take $\frac1{p_4}=\frac 12-\frac{r+\beta-1}3$ when $0<r+\beta-1<\frac
32$ and any $p_4 $ finite
when $r+\beta-1\ge \frac 32$ according to \eqref{Sob-32}); then in that case
we set $\frac 1{p_3}=\frac 12-\frac 1{p_4}\in (2,\infty)$ and use that
\eqref{eq-11} holds for any finite $p_3$ according to \eqref{Sob-32},
since $\beta+2\alpha=\frac 52-\beta\ge \frac 32$.

Hence we have obtained that
\[
\|B(u,v)\|_{V^{r-\beta}}\le C \|v\|_{V^r}  \|v\|_{V^{r+\beta}}.
\]
This implies
\[
\frac{dv}{dt} = -B(u,v)  -\nu\Lambda^{2\beta}v +\Pi \theta e_3 \in L^2(0,T;V^{r-\beta}) .
\]
By Lemma \ref{lemmaC0} we deduce that $v \in C([0,T];V^r)$.

As far as the  continuity in time for $\theta$ is concerned, we have 
that $\theta$ satisfies a transport equation
\[
\partial_t \theta +u \cdot \nabla \theta=0
\]
where the velocity is given and in particular
$u \in C([0,T];V^{r+2\alpha})$ with $r+2\alpha>\frac 52$ (since,  by assumption, $r>2\beta=\frac 52-2\alpha$).
\cite{kp86} considers  this equation
in $\mathbb R^2$; but a straightforward modification of Lemma 4.4
 of \cite{kp86} allows to prove  in the three dimensional case
that given $u \in C([0,T];V^\rho)$ with $\rho>\frac 52$ and $\theta_0 \in H^k$ with $0\le k <[\rho]$, then there exists a unique solution 
$\theta \in C([0,T];H^k)$. Taking $\rho=r+2\alpha$ and $k=r-\beta$, we get the continuity result for $\theta$.
\hfill $\Box$

Now, this regularity is enough to get uniqueness.

\begin{thm}[Uniqueness]\label{uniq} 
We are given  parameters $\alpha$ and $\beta$  with $\frac 12 <\beta<\frac 54$ and 
\[
\alpha+\beta= \frac 54 . 
\]
Let 
\[
r>\max\left(2\beta,\beta+1\right).
\] 
Then, the solutions given in Theorem \ref{Tr} are unique.
\end{thm}
\pf
Let $(v_1,\theta_1)$ and $(v_2,\theta_2)$ be two solutions
 given by Theorem \ref{Tr}.
We define $V=v_1-v_2$, $U=u_1-u_2$ and $\Phi=\theta_1-\theta_2$.
Using the bilinearity we have that they satisfy
\[
\begin{cases}
\partial_t V +\nu \Lambda^{2\beta} V +B(u_1,V)
+B(U,v_2)=\Pi \Phi e_3\\
\partial_t \Phi+U\cdot \nabla \theta_1+u_2\cdot\nabla \Phi=0
\end{cases}
\]
As before, using \eqref{tril0} we get
\[
\begin{split}
\frac 12 \frac{d}{dt}\|V(t)\|_{V^0}^2 &+ \nu \|V(t)\|_{V^{\beta}}^2
\\
&=-\langle B(u_1(t),V(t)),V(t)\rangle
-\langle B(U(t),v_2(t)),V(t)\rangle +\langle \Phi(t) e_3,V(t)\rangle\\
&\le-\langle B(U(t),v_2(t)),V(t)\rangle+\|\Phi(t)\|_{H^0} \|V(t)\|_{V^0}.
\end{split}
\]
And similarly, using \eqref{tril1}
\begin{equation*}
\begin{split}
\frac 12 \frac{d}{dt}\|\Phi(t)\|_{H^{0}}^2 
&=
-\langle (U(t)\cdot \nabla \theta_1(t)),\Phi(t)\rangle
-\langle u_2(t)\cdot \nabla \Phi(t),\Phi(t)\rangle\\
&=
-\langle (U(t)\cdot \nabla \theta_1(t)),\Phi(t)\rangle .
\end{split}
\end{equation*}
Let us estimate the terms on the right hand side 
of each of the relationships above. 
For the velocity equation, we proceed  as usual 
by means of H\"older and Sobolev inequalities
with $\frac 1{p_2}=\frac 12-\frac{2\beta-1}3\in (0,\frac 12)$ and
$\frac 1{p_1}=\frac 12-\frac 1{p_2}\equiv \frac 12-\frac{2\alpha}3$:
\begin{equation*}
\begin{split}
\left | \langle B(U,v_2),V\rangle\right | 
&\le \|(U\cdot \nabla)v_2\|_{L_2}\|V\|_{L_2}\\
& \le \|U\|_{L_{p_1}} \|\nabla v_2\|_{L_{p_2}}\|V\|_{V^0} \\
&\le C \|U\|_{V^{2\alpha}} \|v_2\|_{V^{2\beta}}\|V\|_{V^0} \\
&=C \|V\|_{V^0} \|v_2\|_{V^{2\beta}}\|V\|_{V^0}\\
&\le C \|V\|_{V^{\beta}} \|v_2\|_{V^{2\beta}}\|V\|_{V^0}\\
&\le \frac{\nu}{4} \|V\|_{V^{\beta}}^2+C_\nu \|v_2\|^2_{V^{2\beta}}\|V\|^2_{V^0}
\end{split}
\end{equation*}

Similarly, for the temperature equation:
\[
\left |\langle U\cdot \nabla \theta_1,\Phi \rangle \right |
\le \|U\cdot \nabla\theta_1\|_{L_2}\|\Phi\|_{L_2}
 \le \|U\|_{L_{p_3}} \|\nabla \theta_1\|_{L_{p_4}}\|\Phi\|_{H^0}
\]
with $\frac 1{p_3}+\frac 1{p_4}=\frac 12$. Now we choose $p_3$ and $p_4$.
When $1<\beta<\frac 54$ we set $\frac 1{p_3}=\frac 12-\frac{\frac 52-\beta}3$
and $\frac 1{p_4}=\frac
12-\frac{\beta-1}3$ so to get 
$\|U\|_{L_{p_3}}\le C \|U\|_{V^{\frac 52-\beta}}$ and
$\|\nabla \theta_1\|_{L_{p_4}}\le C \|\theta_1\|_{H^{\beta}}$; in
addition we use that $H^{r-\beta}\subseteq H^\beta$ when $r\ge 2\beta$. Therefore
\begin{equation}\label{stima-u-theta}
\|U\|_{L_{p_3}} \|\nabla \theta_1\|_{L_{p_4}}\le  C\|U\|_{V^{\frac 52-\beta}}\|\theta_1\|_{H^{r-\beta}}
=C \|V\|_{V^\beta}\|\theta_1\|_{H^{r-\beta}}.
\end{equation}
On the other hand, when $\beta \le 1$, according to \eqref{Sob-32} we have
$\|U\|_{L_{p_3}}\le C \|U\|_{V^{\frac 52-\beta}}$
for any finite $p_3$; hence we first choose $p_4>2$ such that
$\|\nabla \theta_1\|_{L_{p_4}}\le C\|\nabla \theta_1\|_{H^{r-\beta-1}}\le C \|\theta_1\|_{H^{r-\beta}}$;
this can be done as soon as $r-\beta-1>0$, i.e. $r>\beta+1$
(as in \eqref{nuovo-argomento}). 
Then we set $\frac 1{p_3}=\frac 12-\frac 1{p_4}$. Again we have obtained \eqref{stima-u-theta}.  
\\
Thus
\[
\left |\langle U\cdot \nabla \theta_1,\Phi \rangle \right |
\le C \|V\|_{V^{\beta}} \|\theta_1\|_{H^{r-\beta}}\|\Phi\|_{H^0}
\le \frac{\nu}{4} \|V\|_{V^{\beta}}^2+C_\nu  \|\theta_1\|^2_{H^{r-\beta}}\|\Phi\|^2_{H^0}.
\]

Summing up, we have obtained 
\begin{multline*}
\frac{d}{dt}\|V(t)\|_{V^0}^2+\nu \|V(t)\|_{V^{\beta}}^2
+ \frac{d}{dt}\|\Phi(t)\|_{H^{0}}^2 \\
\le
C \|v_2(t)\|^2_{V^{2\beta}}\|V(t)\|^2_{V^0}
+ \|\theta_1(t)\|^2_{H^{r-\beta}}\|\Phi(t)\|^2_{H^0}+
\|\Phi(t)\|_{H^{0}}^2+\|V(t)\|_{V^0}^2.
\end{multline*}
If we define  $Z(t)=\|V(t)\|_{V^0}^2+\|\Phi(t)\|^2_{H^0}$, we have $Z(0)=0$ and
\[
Z^\prime(t)\le C(\|v_2(t)\|^2_{V^{2\beta}}+ \|\theta_1(t)\|^2_{H^{r-\beta}}+ 1)Z(t).
\]
By Gronwall lemma we get $Z(t)=0$ for all $t$, and this completes the proof.
\hfill $\Box$

\begin{thm}[Continuous dependence on the initial data] \label{cc}
We are given  parameters $\alpha$ and $\beta$  with $\frac 12 <\beta<\frac 54$ and 
\[
\alpha+\beta= \frac 54 . 
\]
Let
\[
r> \beta+2.
\]
Then, given any 
initial conditions $v_{1,0},v_{2,0} \in V^r$ and $\theta_{1,0},\theta_{2,0} \in H^{r-\beta}$
we have
\begin{multline}
\|v_1-v_2\|_{L^\infty(0,T;V^{r-1})}+\|v_1-v_2\|_{L^2(0,T;V^{r-1+\beta})}
+\|\theta_1-\theta_2\|_{L^\infty(0,T;H^{r-\beta-1})}
\\\le C \left(\|v_{1,0}-v_{2,0}\|_{V^{r-1}}+ \|\theta_{1,0}-\theta_{2,0}\|_{H^{r-\beta-1}}\right) 
\end{multline}
where the constant $C$ depends on $T$, 
$\|\theta_1\|_{L^\infty(0,T;H^{r-\beta})}$, $\|v_i\|_{L^2(0,T;V^{r+\beta-1})}$ and 
$\|v_i\|_{L^\infty(0,T;V^{r})}$.
\end{thm}
\pf
We begin by pointing out that, under the assumption $\frac
12<\beta<\frac 54$ the condition $r> \beta+2$ implies also
$r>\max(2\beta,\beta+1,2-\beta)$ and therefore the assumptions of
Theorem \ref{Tr} and 
Lemma \ref{lBdiff}, \ref{t-diff} and \ref{t-diff2} are fulfilled.

Using the same setting as in the proof of  Theorem \ref{uniq}, we get
\[
\begin{split}
\frac 12 \frac{d}{dt}&\|V(t)\|_{V^{r-1}}^2 + \nu \|V(t)\|_{V^{r-1+\beta}}^2
=-\langle B(\Lambda^{-2\alpha}v_1(t),V(t)),\Lambda^{2r-2} V(t)\rangle\\
&
\;\;-\langle B(\Lambda^{-2\alpha} V(t),v_2(t)),\Lambda^{2r-2} V(t)\rangle 
+\langle \Lambda^{r-\beta-1}\Phi(t) e_3,\Lambda^{r-1+\beta}V(t)\rangle.
\end{split}
\]
We estimate the first two terms of r.h.s. by means of Lemma \ref{lBdiff}
\[
\begin{split}
&|\langle B(\Lambda^{-2\alpha}v_1(t),V(t)),\Lambda^{2r-2} V(t)\rangle|
+|\langle B(\Lambda^{-2\alpha} V(t),v_2(t)),\Lambda^{2r-2} V(t)\rangle|
\\
&\quad\le C\big( \|v_1\|_{V^{r}} \|V\|_{V^{r-1}}
               +\|v_1\|_{V^{r+\beta-1}}\|V\|_{V^{r+\beta-1}}\big)\|V\|_{V^{r-1}}\\
 & \qquad +C\|V\|_{V^{r-1}} \|v_2\|_{V^{r+\beta-1}} \|V\|_{V^{r+\beta-1}}.
          %+\|\Phi\|_{H^{r-\beta-1}}\|V\|_{V^{r-1+\beta}}
\end{split}
\]
Using Young inequality, we get
\begin{multline}
\frac 12 \frac{d}{dt}\|V(t)\|_{V^{r-1}}^2 + \nu \|V(t)\|_{V^{r-1+\beta}}^2
\le
\frac \nu 2 \|V(t)\|_{V^{r+\beta-1}}^2 + C_\nu  \|\Phi(t)\|_{H^{r-\beta-1}}^2 \\
+ C_\nu (\|v_1(t)\|_{V^{r}}+\|v_1(t)\|_{V^{r+\beta-1}}^2
+\|v_2(t)\|_{V^{r+\beta-1}}^2)\|V(t)\|_{V^{r-1}}^2.
\end{multline}

Similarly,
for the temperature difference; we use Lemma \ref{t-diff} and \ref{t-diff2} and Young 
 inequality
\[\begin{split}
\frac 12 \frac{d}{dt}\|\Phi(t)\|_{H^{r-\beta-1}}^2 
&=
-\langle U(t)\cdot \nabla \theta_1(t),\Lambda^{2r-2\beta-2}\Phi(t)\rangle
-\langle u_2(t)\cdot \nabla \Phi(t),\Lambda^{2r-2\beta-2}\Phi(t)\rangle
\\
&\le 
C \|V(t)\|_{V^{r-1}}\|\theta_1(t)\|_{H^{r-\beta}}\|\Phi(t)\|_{H^{r-\beta-1}}
+C \|v_2(t)\|_{V^{r}}\|\Phi(t)\|_{H^{r-\beta-1}}^2
\\
&\le C \|V(t)\|_{V^{r-1}}^2
     + C \big( \|\theta_1(t)\|_{H^{r-\beta}}^2
+ \|v_2(t)\|_{V^r}\big)\|\Phi(t)\|_{H^{r-\beta-1}}^2.
\end{split}
\]

Finally, we consider the sum $\|V(t)\|_{V^{r-1}}^2 +\|\Phi(t)\|_{H^{r-\beta-1}}^2:=W(t)$ and define
$a(t)=1+\|\theta_1(t)\|_{H^{r-\beta}}^2
 +\|v_1(t)\|^2_{V^{r+\beta-1}}+\|v_2(t)\|^2_{V^{r+\beta-1}}+
 \|v_1(t)\|_{V^{r}}+\|v_2(t)\|_{V^{r}}$;
we  have $a \in L^1(0,T)$ and 
\begin{equation}\label{stimo-le-diff}
W^\prime(t)
    + \nu \|V(t)\|_{V^{r+\beta-1}}^2
\le C a(t) W(t).
\end{equation}
Gronwall lemma applied to
\[
W^\prime(t)\le C a(t) W(t)
\] 
gives 
\[
\sup_{0\le t \le T}W(t)\le W(0)e^{C\int_0^T a(t)\ dt} .
\]
Integrating in time \eqref{stimo-le-diff} and using the latter result 
we get the estimate
for $\int_0^T \|V(t)\|_{H^{r+\beta-1}}^2dt$.
This concludes the proof.
\hfill $\Box$

\section{Auxiliary results}
In this section we prove the  lemma used in the proofs of the previous section.
\begin{lemma}\label{lv1}
Let $\frac 12<\beta < \frac 54$ and $\alpha+\beta= \frac 54$.
Then, for any $r>0$ there exists a constant $C>0$ such that
\[
|\langle B(\Lambda^{-2\alpha} v,v),\Lambda^{2r}v\rangle|
\le C
 \|v\|_{V^{2\beta}} \|v\|_{V^{r+\beta}} \|v\|_{V^{r}}.
\]
\end{lemma}
\pf
Set  $u = \Lambda^{-2\alpha}v$.
First 
\[
\begin{split}
\langle B(u, v),\Lambda^{2r}v\rangle
&=\langle \Lambda^r\Big((u \cdot \nabla) v\Big),\Lambda^{r}v\rangle
\\
&=
\langle [\Lambda^r,u ] \cdot\nabla v,\Lambda^{r}v\rangle \qquad
  \text{ by } \eqref{tril-con-comm2}\\
&\le \|[\Lambda^r,u ] \cdot\nabla v\|_{L_2} \|\Lambda^{r}v\|_{L_2}.
\end{split}
\]
Then, we use the Commutator Lemma \ref{comm} with $p=2$ and
\[
\tfrac 1{p_1}=\tfrac 12-\tfrac{\beta-1+2\alpha}3,\qquad
\tfrac 1{p_2}=\tfrac 12-\tfrac 1{p_1}\equiv \tfrac 12-\tfrac \beta 3 \in (\tfrac 1{12},\tfrac 13)
\]
\[
\begin{cases}
\frac 1{p_3}=\frac 12-\frac{\beta+2\alpha}3 \in (0,\frac1{12}), \qquad
\frac 1{p_4}=\frac 12-\frac 1{p_3}\equiv \frac 12-\frac {\beta-1}3
    & \text{ if } 1<\beta<\frac 54\\
\frac 1{p_3}=\frac 12-\frac{2\alpha}3 \in (0,\frac 13], \frac 1{p_4}=\frac 12-\frac 1{p_3}\equiv\frac 12-\frac{2\beta-1}3
   & \text{ if } \frac 12<\beta\le 1
\end{cases}
\]
so to get
\[\begin{split}
\|[\Lambda^r,u ]\cdot \nabla v\|_{L_2} &\le C\left(\|\Lambda u\|_{L_{p_1}} \|\Lambda ^r v\|_{L_{p_2}}+
     \|\Lambda^r u\|_{L_{p_3}}\|\Lambda v\|_{L_{p_4}}\right)\\
     & =C\left(\|\Lambda^{1-2\alpha} v\|_{L_{p_1}} \|\Lambda ^r v\|_{L_{p_2}}+
     \|\Lambda^{r-2\alpha} v\|_{L_{p_3}}\|\Lambda v\|_{L_{p_4}}\right).
\end{split}\]
Then, to conclude our estimate we use  the Sobolev embedding inequalities
\begin{equation}\label{stima:v-1-2}
\|\Lambda^{1-2\alpha} v\|_{L_{p_1}}\le C \|v\|_{V^\beta} %\le C \|v\|_{V^{2\beta} }
\qquad
\|\Lambda ^r v\|_{L_{p_2}}\le C \|v\|_{V^{r+\beta}}
\end{equation}
and for $1<\beta<\frac 54$\[
\|\Lambda^{r-2\alpha} v\|_{L_{p_3}}\le C \|v\|_{V^{r+\beta}} \qquad
\|\Lambda  v\|_{L_{p_4}}\le C \|v\|_{V^\beta} \le C \|v\|_{V^{2\beta}},
\]
whereas for $\frac 12<\beta\le 1$ (i.e. $\frac 14\le \alpha<\frac 34$) 
\[
\|\Lambda^{r-2\alpha} v\|_{L_{p_3}}\le C \|v\|_{V^r}\le C\|v\|_{V^{r+\beta}}
\qquad
 \|\Lambda  v\|_{L_{p_4}}\le C \|v\|_{V^{2\beta}}.
\]
\hfill $\Box$

\begin{lemma}\label{lBdiff}
Let $\frac 12<\beta < \frac 54$ and $\alpha+\beta= \frac 54$.
If
\[
r>\max(2\beta, 2-\beta) ,
\]
then there exists a constant $C>0$ such that
\[
|\langle B(\Lambda^{-2\alpha} w,v),\Lambda^{2r-2}v\rangle|
\le C
 (\|w\|_{V^{r}} \|v\|_{V^{r-1}}+\|w\|_{V^{r+\beta-1}}\|v\|_{V^{r+\beta-1}})
  \|v\|_{V^{r-1}}
\]
and
\[
|\langle B(\Lambda^{-2\alpha} v,w),\Lambda^{2r-2}v\rangle|
\le C
 \|v\|_{V^{r-1}}\|w\|_{V^{r+\beta-1}} \|v\|_{V^{r+\beta-1}}.
\]
\end{lemma}
\pf
First, notice that we also have $r> 1$.\\
To prove the first inequality, we 
use the Commutator Lemma \ref{comm} with $p=p_2=2$, $p_1=\infty$ and suitables $p_3, p_4$ 
to get
\[\begin{split}
|\langle B(\Lambda^{-2\alpha}& w,v),\Lambda^{2r-2}v\rangle|
\\&=|\langle \Lambda^{r-1}\Big((\Lambda^{-2\alpha} w\cdot\nabla) v\Big), \Lambda^{r-1}v\rangle|
\\
&=
|\langle [\Lambda^{r-1},\Lambda^{-2\alpha} w]\cdot \nabla  v, \Lambda^{r-1}v\rangle|
\qquad \text{ by } \eqref{tril-con-comm2}\\
&\le C
(\|\Lambda^{1-2\alpha} w\|_{L_\infty}\|\Lambda^{r-1}v\|_{L_2}
+\|\Lambda^{r-1-2\alpha}w\|_{L_{p_3}}\|\Lambda v\|_{L_{p_4}})\|\Lambda ^{r-1}v\|_{L_2}
\end{split}
\]
We estimate the first four terms in the latter line.
When $1<\beta<\frac 54$ we choose
\begin{equation}\label{p34-1}
\tfrac 1{p_3}=\tfrac 12 -\tfrac {\beta+2\alpha}3 \in (0,\tfrac 1{12}),\qquad
\tfrac 1{p_4}=\tfrac 12 -\tfrac 1{p_3}\equiv \tfrac 12-\tfrac{\beta-1}3,
\end{equation}
whereas when $\frac 12<\beta \le 1$ we choose 
\begin{equation}\label{p34-2}
\begin{cases}
\frac 1{p_4}=\frac 12-\frac{r+\beta-2}3 & \text{ if } 2<r+\beta<\frac 72
\\
\text{ any } p_4\in (2,\infty) & \text{ if } r+\beta\ge\frac 72
\end{cases}
\end{equation}
and $\frac 1{p_3}=\frac 12-\frac 1{p_4} \in (0,\frac 12)$.

Then we use the following Sobolev embedding inequalities:
\[
\|\Lambda^{1-2\alpha} w\|_{L_\infty}\le C \|w\|_{V^{r}}
\]
since $r>\frac 52-2\alpha=2\beta$. Moreover, for $1<\beta<\frac 54$, according to \eqref{p34-1} we have
\[
\|\Lambda^{r-1-2\alpha}w\|_{L_{p_3}}\le C \|w\|_{V^{r+\beta-1}}
\qquad
\|\Lambda v\|_{L_{p_4}} \le C \|v\|_{V^\beta}\
\]
and 
\[
\|v\|_{V^\beta}\le C \|v\|_{V^{r-1+\beta}}
\]
since $r-1>0$.  On the other hand, 
for $\frac 12<\beta\le 1$ according to \eqref{p34-2} there exists $p_4\in (2,\infty)$ such that
\[
\|\Lambda v\|_{L_{p_4}} \le C \|v\|_{V^{r+\beta-1}};
\]
then we set $\frac 1{p_3}=\frac 12-\frac 1{p_4}\in (0,\frac 12)$ and recall that 
\begin{equation}\label{p34-3}
\|\Lambda^{r-1-2\alpha}w\|_{L_{p_3}}=\|\Lambda^{r-\frac 72+2\beta}w\|_{L_{p_3}}
\le C \|w\|_{V^{r+2\beta-2}} 
\end{equation}
for any finite  $p_3$. Since $\|w\|_{V^{r+2\beta-2}} \le C \|w\|_{V^{r+\beta-1}}$ when $\beta\le 1$,
this concludes the first inequality of the statement of this Lemma.

For the second inequality, we use Lemma \ref{prod} with $p=2$:
\begin{multline*}
|\langle B(\Lambda^{-2\alpha} v,w),\Lambda^{2r-2}v\rangle|
=|\langle \Lambda^{r-1-\beta}\Big((\Lambda^{-2\alpha} v\cdot \nabla )w\Big), 
  \Lambda^{r+\beta-1}v\rangle|
\\
\le C
\Big(\|\Lambda^{-2\alpha}v\|_{L_{p_1}}\|\Lambda^{r-\beta}w\|_{L_{p_2}}+
  \|\Lambda^{r-1-\beta-2\alpha}v\|_{L_{p_3}} \|\Lambda w\|_{L_{p_4}}\Big)
\|v\|_{V^{r+\beta-1}}.
\end{multline*}
Now we choose 
$\frac 1{p_1}=\frac 12-\frac{2\alpha}3\equiv\frac {2\beta-1}3$ and 
$\frac 1{p_2}=\frac 12-\frac{2\beta-1}3 \in (0,\frac 12)$ 
since $\frac 12<\beta<\frac 54$; then, 
by means of  Sobolev embedding inequalities
\[
\|\Lambda^{-2\alpha}v\|_{L_{p_1}}\le C\|v\|_{V^0}\le C \|v\|_{V^{r-1}} \;\text{ for any } r\ge 1
\]
\[
\|\Lambda^{r-\beta}w\|_{L_{p_2}}\le C \|w\|_{V^{r+\beta-1}}
\]
Moreover, for $1<\beta<\frac 54$ we choose
$\frac 1{p_3}=\frac 12-\frac{\beta+2\alpha}3\equiv \frac{\beta-1}{3}\in (0,\frac 1{12})$ and 
$\frac 1{p_4}=\frac 12-\frac{\beta-1}{3}$; therefore
by means of Sobolev embedding theorems we get
\[
\|\Lambda^{r-1-\beta-2\alpha}v\|_{L_{p_3}} \le
C \|v\|_{V^{r-1}}
\]
\[
\|\Lambda w\|_{L_{p_4}} \le
C \|w\|_{V^{\beta}}\le C \|w\|_{V^{r-1+\beta}} \text{ for any } r \ge 1
\]
On the other side, for $\frac 12<\beta\le 1$, 
we choose $p_4\in (2,\infty)$ as in \eqref{p34-2} so to get
\[
\|\Lambda w\|_{L_{p_4}} \le C \|w\|_{V^{r-1+\beta}}
\]
and we set $\frac 1{p_3}=\frac 12-\frac 1{p_4}$ to get, as in \eqref{p34-3},
\[
\|\Lambda^{r-1-\beta-2\alpha}w\|_{L_{p_3}}=\|\Lambda^{r-\frac 72+\beta}w\|_{L_{p_3}}
\le C \|w\|_{V^{r+\beta-2}} . 
\]
Using that 
\[ \|v\|_{V^{r+\beta-2}}\le C \|v\|_{V^{r-1}}
\]
since $\beta\le 1$, we  conclude the second inequality in the statement. \hfill $\Box$

\begin{lemma}\label{t-diff}
Let $\frac 12<\beta < \frac 54$ and $\alpha+\beta= \frac 54$.
If 
\[
r>\max(2\beta, \beta+1)
\]
then there exists a constant $C>0$ such that
\[
|\langle \Lambda^{-2\alpha} v \cdot \nabla \theta,\Lambda^{2r-2\beta-2}\phi\rangle|
\le C
  \| v\|_{V^{r-1}} \|\theta\|_{H^{r-\beta}} \|\phi\|_{H^{r-\beta-1}}.
\]
\end{lemma}
\pf
We use Lemma  \ref{prod} with  
$p=p_2=2$ and $p_1=\infty$:
\[
\begin{split}
|\langle \Lambda^{-2\alpha} &v \cdot \nabla \theta,\Lambda^{2r-2\beta-2}\phi\rangle|
\\&=
|\langle \Lambda^{r-\beta-1}(\Lambda^{-2\alpha} v \cdot \nabla \theta),
  \Lambda^{r-\beta-1}\phi\rangle|
\\
&\le
\|\Lambda^{r-\beta-1}(\Lambda^{-2\alpha} v \cdot \nabla \theta)\|_{L_2}
\|\phi\|_{H^{r-\beta-1}}
\\
&\le
C (\|\Lambda^{-2\alpha} v\|_{L_\infty} \|\Lambda^{r-\beta}\theta\|_{L_{2}}+
\|\Lambda^{r-\beta-1-2\alpha} v\|_{L_{p_3}} \|\Lambda\theta\|_{L_{p_4}})
\|\phi\|_{H^{r-\beta-1}}
\end{split}
\]
We estimate the first four terms in the latter line.
Since $r>2\beta$, i.e. $r-1-(2\beta-\frac 52)>\frac32$, we have
\[
\|\Lambda^{-2\alpha} v\|_{L_\infty}= \|\Lambda^{2\beta - \frac 52} v\|_{L_\infty}
\le C \|v\|_{V^{r-1}} .
\]
For $1<\beta<\frac 54$ we set
$\frac 1{p_3}=\frac12-\frac{\beta+2\alpha}3\in (0,\frac1{12})$ 
and $\frac 1{p_4}=\frac 12-\frac 1{p_3}\equiv 
\frac 12-\frac{\beta-1}3$
so to get
\[
\|\Lambda^{r-\beta-1-2\alpha} v\|_{L_{p_3}}
\le C \|v\|_{V^{r-1}}
\]
\[
\|\Lambda\theta\|_{L_{p_4}}\le C \|\theta\|_{H^\beta} \le C  \|\theta\|_{H^{r-\beta}} \; \text{ when }  r\ge 2\beta
\]
On the other side, when $\frac 12<\beta\le 1$ we have $\beta+2\alpha\ge \frac 32$; 
hence, according to \eqref{Sob-32}
\[
\|\Lambda^{r-\beta-1-2\alpha} v\|_{L_{p_3}}\le C \|v\|_{V^{r-1}}
\qquad \text{ for any finite } p_3 .
\]
Therefore we set $\frac 1{p_3}=\frac 12-\frac 1{p_4}$ with $p_4\in (2,\infty)$ chosen arbitrarily 
when $r\ge \beta+\frac 52$ and $\frac 1{p_4}=\frac 12-\frac{r-\beta-1}3$ when 
$\beta+1< r< \beta+\frac 52$; in this way we get
\[
\|\Lambda\theta\|_{L_{p_4}}\le C \|\theta\|_{H^{r-\beta}} .
\]
This concludes the proof.
\hfill $\Box$

\begin{lemma}\label{t-diff2}
Let $\frac 12<\beta < \frac 54$ and $\alpha+\beta= \frac 54$.
If 
\[
 r> \beta+2
\]
then there exists a constant $C>0$ such that
\[
|\langle\Lambda^{-2\alpha} v\cdot\nabla \theta,\Lambda^{2r-2\beta-2 }\theta\rangle|
\le C
 \|v\|_{V^{r}}\|\theta\|_{H^{r-\beta-1}}^2.
\]
\end{lemma}
\pf
We use the Commutator 
Lemma \ref{comm} with $p=p_2=2$,
$p_1=\infty$:
\[\begin{split}
|\langle\Lambda^{-2\alpha} &v\cdot\nabla \theta,\Lambda^{2r-2\beta-2 }\theta\rangle|
\\&
=|\langle\Lambda^{r-\beta-1}(\Lambda^{-2\alpha} v\cdot\nabla \theta),
 \Lambda^{r-\beta-1}\theta\rangle|
\\
&=
|\langle[\Lambda^{r-\beta-1}, \Lambda^{-2\alpha} v]\cdot\nabla \theta,
 \Lambda^{r-\beta-1}\theta\rangle| \qquad \text{ by } \eqref{tril-con-comm1}
\\
&\le C(\|\Lambda^{1-2\alpha}v\|_{L_\infty} \|\Lambda^{r-\beta-1}\theta\|_{L_2}+
\|\Lambda^{r-\beta-1-2\alpha} v\|_{L_{p_3}} \|\Lambda\theta\|_{L_{p_4}})
\|\theta\|_{H^{r-\beta-1}}
\end{split}
\]
We estimate the first four terms in the latter line.
Since $r>\beta+2>2\beta-1$ we have
\[
\|\Lambda^{1-2\alpha}v\|_{L_\infty}=\|\Lambda^{2\beta-\frac 52}v\|_{L_\infty}
\le C \|v\|_{V^r} .
\]
Moreover  we have
\[
\|\Lambda^{r-\beta-1-2\alpha} v\|_{L_{p_3}}\le 
C\|\Lambda^{r-\frac 94} v\|_{L_{p_3}}
\]
and according to \eqref{Sob-32}
\[
\|\Lambda^{r-\frac 94} v\|_{L_{p_3}}\le C \|v\|_{V^r}
\]
for any finite $p_3$.
Hence we set $\frac 1{p_3}=\frac 12-\frac 1{p_4}$ with $p_4 \in (2,\infty)$ chosen arbitrarily when 
$r\ge \beta+\frac 72$  and $\frac 1{p_4}=\frac 12-\frac{r-\beta-2}3$ when 
$\beta+2<r<\beta+\frac 72$ in order to have the Sobolev inequality
\[
\|\Lambda\theta\|_{L_{p_4}}\le C \|\theta\|_{H^{r-\beta-1}}.
\]
\hfill $\Box$

\begin{lemma}\label{ltt}
Let $\alpha+\beta= \frac 54$ with $\frac 12<\beta < \frac 54$  and  $r>\beta+1$.
Then, there exist  $q_3, q_4 >2$ with $\frac 1{q_3}+\frac 1{q_4}\le \frac 12$
and  $q >2$,  $a \in (0,1)$, $C>0$ such that
\[
\|\Lambda^{r-\beta-2\alpha}v\|_{L_{q_3}}
\|\Lambda \theta\|_{L_{q_4}}
\|\theta\|_{H^{r-\beta}}
\le C
\|v\|_{V^{2\beta}}^{a} \|v\|_{V^{r+\beta}}^{1-a} 
\|\theta\|_{L_q}^{1-a}\|\theta\|_{H^{r-\beta}}^{1+a}.
\]
\end{lemma}
\pf
We use Sobolev embedding theorem,  interpolation theorem
 and the Gagliardo-Nirenberg inequality; 
then for some $a \in (0,1)$  and $q\ge 2$ to be defined later on we look for  
\[\begin{cases}
\|\Lambda^{r-\beta-2\alpha}v\|_{L_{q_3}}
\le C\|\Lambda ^{r+\beta-ra +\beta a}v\|_{L_2}\le C\|v\|_{V^{2\beta}}^{a} \|v\|_{V^{r+\beta}}^{1-a}
& \text{ for } \frac 1{q_3}=\frac 12-\frac{2\beta+\beta a-r a+2\alpha}3\\
\|\Lambda \theta\|_{L_{q_4}}
\le C \|\theta\|_{L_q}^{1-a} \|\Lambda^{r-\beta}\theta\|^{a}_{L_2}
& \text{ for } \frac 1{q_4}=\frac 13+(\frac 12 -\frac{r-\beta}3)a+\frac {1-a}q
\end{cases}
\]
under the conditions
\[
\begin{cases}
r+\beta-ra+\beta a\ge r-\beta-2\alpha
\\
\frac 1{q_3}+\frac 1{q_4}\le \frac 12
\\
\frac{1}{r-\beta}< a <1
\end{cases}
\]
equivalent to (since $r>\beta$ by assumption)
\begin{equation}\label{sistema-A}
\begin{cases}
a\le \frac 5{2(r-\beta)}\\
\frac a2+\frac{1-a}q\le \frac 12\\
\frac{1}{r-\beta}< a <1
\end{cases}
\end{equation}
The second equation is satisfied for some $q$ (big enough) when $0<a<1$;
therefore we choose  $a \in (0,1)$ such that
\begin{equation}\label{condizione-a}
\frac 1{r-\beta}< a < 
\min\Big( 1, \frac 5{2(r-\beta)}\Big).
\end{equation}
This double condition has solutions since $r-\beta>1$.
\hfill $\Box$
\\[2mm]
{\bf Acknowledgements}
The research of Hakima Bessaih was supported by the NSF grants DMS-1416689 
and DMS-1418838. 
Part of this research started while Hakima Bessaih was visiting the 
Department of Mathematics of the University of Pavia and was partially 
supported by  the GNAMPA-INdAM Project 2014 ''Regolarit\`a e dissipazione in 
fluidodinamica'';  she would like to thank the hospitality of the Department.  

We are very grateful to the anonymous referee;
his/her careful reading and suggestions helped to improve greatly  
the final result of the paper.


\begin{thebibliography}{99}

   \bibitem{AH07}
H. Abidi, T. Hmidi:
On the global well-posedness for Boussinesq system, 
{\it J. Differential Equations} 233 (2007), no. 1, 199-220

   \bibitem{AHK11}
H. Abidi, T. Hmidi, S. Keraani:
On the global regularity of axisymmetric Navier-Stokes-Boussinesq system,
{\it Discrete Contin. Dyn. Syst.} 29 (2011), no. 3, 737-756

   \bibitem{BBF}
D. Barbato, H. Bessaih, B. Ferrario:
On a stochastic Leray-$\alpha$ model of Euler equations,
{\it Stochastic Process. Appl.} {\bf 124} (2014), no. 1,  199-219

   \bibitem{BMR}
D. Barbato, F. Morandin, M. Romito:
Global regularity for a slightly supercritical 
hyperdissipative Navier-Stokes system,
{\it Analysis and PDE} {\bf 7} (2014), no. 8, 2009-2027

\bibitem{BG}
H. Br\'ezis, T. Gallouet:
 Nonlinear Schr\"odinger evolution equations,
 {\it Nonlinear Anal.} 4 (1980), no. 4, 677-681
 
    \bibitem{BW}
H. Br\'ezis, S. Wainger:
A note on limiting cases of Sobolev embeddings and convolution inequalities,
{\it Comm. Partial Differential Equations} 5 (1980), no. 7, 773-789

   \bibitem{C06}
D. Chae:
Global regularity for the 2D Boussinesq equations with partial viscosity terms,
{\it Adv. Math}. 203 (2006), no. 2, 497-513
      
   \bibitem{DP08}
R. Danchin, M. Paicu:
Existence and uniqueness results for the Boussinesq system with data in Lorentz spaces,
{\it Phys. D} 237 (2008), no. 10-12, 1444-1460
   
   \bibitem{DP09}
R. Danchin, M. Paicu: 
Global well-posedness issues for the inviscid Boussinesq system with Yudovich's type data,
{\it Comm. Math. Phys.} 290 (2009), no. 1, 1-14

   \bibitem{GF12}            
J. Geng, J. Fan:
A note on regularity criterion for the 3D Boussinesq system with 
zero thermal conductivity,
{\it Appl. Math. Lett.} 25 (2012), no.~1, 63-66

   \bibitem{HK07}
T. Hmidi, Keraani:
On the global well-posedness of the two-dimensional Boussinesq system 
with a zero diffusivity,
{\it  Adv. Differential Equations} 12 (2007), no. 4, 461-480

    \bibitem{HL05}
T. Y. Hou, C. Li:
Global well-posedness of the viscous Boussinesq equations,
{\it Discrete Contin. Dyn. Syst.} 12 (2005), no. 1, 1-12 
   
   \bibitem{kp86}
T. Kato, G. Ponce:
Well posedness of the Euler and Navier-Stokes equations in the Lebesgues spaces $L^p_s(\mathbb R^2)$,
{\it Rev. Mat. Iberoam. } 2 (1986), no. 1-2, 73-88 

   \bibitem{kp}
T. Kato, G. Ponce:
Commutator estimates and the Euler and Navier-Stokes equations, 
{\it Comm. Pure Appl. Math.} 41  (1988), no. 7, 891-907

  \bibitem{KC}
  M. Kaya, O. \c{C}elebi:
Global attractor for the regularized B\'enard problem,
{\it Applicable Analysis} 93 (2014), no. 9, 1989-2001
  
   \bibitem{kpv}
C. Kenig, G. Ponce, L. Vega:
Well-posedness of the initial value problem for the Korteweg-de-Vries 
equation, 
{\it J. Amer. Math. Soc.} 4 (1991), 323-347

   \bibitem{Lions}
J.-L. Lions: 
{\it Quelques m\'ethodes de r\'esolution des probl\`emes aux limites 
non-lin\'eaires}, Dunod, Gauthier-Villars, Paris 1969

   \bibitem{MS}
J. C. Mattingly,  Ya. G. Sinai:
An elementary proof of the existence and uniqueness theorem for the 
Navier-Stokes equations,
{\it Commun. Contemp. Math.} 1 (1999), no. 4, 497-516 

   \bibitem{Ni}
L. Nirenberg:
On elliptic partial differential equations,
{\it Ann. Scuola Norm. Sup. Pisa}  13 (1959), no. 2, 115-162

   \bibitem{ot}
E. Olson, E. S. Titi:
Viscosity versus vorticity stretching: global well-posedness for a 
family of Navier-Stokes-alpha-like models,
{\it Nonlinear Anal.} 66 (2007), no. 11, 2427-2458

   \bibitem{QDY10}
H. Qiu, Y.  Du, Z. Yao:
A blow-up criterion for 3D Boussinesq equations in Besov spaces,
{\it Nonlinear Anal.} 73 (2010), no. 3, 806-815

  \bibitem{Se}
R. Selmi:  
Global Well-Posedness and Convergence Results for the 3D-Regularized Boussinesq System,
{\it Canad. J. Math.} 64 (2012), no. 6, 1415-1435
  
    \bibitem{St} 
W. A.  Strauss:
On continuity of functions with values in various Banach spaces, 
{\it Pacific J. Math.} 19 (1966), 543-551
 
    \bibitem {temam} 
R. Temam:
{\it  Navier-Stokes equations. Theory and numerical analysis.
Studies in Mathematics and its Applications 2}, 
North-Holland Publishing Co., Amsterdam-New York, 1979.

     \bibitem{temamP} 
R. Temam: 
{\it Navier-Stokes equations and nonlinear functional analysis}, Second edition. 
CBMS-NSF Regional Conference Series in Applied Mathematics, 66. 
Society for Industrial and Applied Mathematics (SIAM), Philadelphia, PA, 1995.

   \bibitem{Kazuo} 
K. Yamazaki: 
On the global regularity of generalized Leray-alpha type models,
{\it Nonlin. Anal.} 75 (2012), 503--515

\bibitem{Ye}
Z. Ye:
A note on global well-posedness of solutions to Boussinesq equations with 
fractional dissipation,
{\it Acta Math. Sci.} Ser. B Engl. Ed. 35 (2015), no. 1, 112-120 

    \bibitem{ZW13}
Xiang Zhaoyin, Yan Wei:
Global regularity of solutions to the Boussinesq equations with fractional diffusion,
{\it Adv. Differential Equations} 18 (2013), no. 11-12, 1105-1128 

\end{thebibliography}
\end{document}